\documentclass[12pt]{article}
\usepackage{graphicx,tikz,color,url}
\usepackage{amsmath,amssymb,latexsym}
\usetikzlibrary{arrows}
\newtheorem{theorem}{Theorem}[section]

\newtheorem{lemma}[theorem]{Lemma}
\newtheorem{proposition}[theorem]{Proposition}

\newcommand{\proof}{\noindent{\bf Proof.\ }}
\newcommand{\qed}{\hfill $\square$ \bigskip}

\newcommand{\gsmb}{\gamma_{\rm SMB}}
\newcommand{\gmb}{\gamma_{\rm MB}}

\newcommand{\gmbt}{\gamma_{\rm MBT}}

\begin{document}
\date{}
\title{Maker-Breaker  total domination number}
\author{ Athira Divakaran$^{a,}$\thanks{Email: \texttt{athiradivakaran91@gmail.com}}
\and
Tijo James$^{b,}$\thanks{Email: \texttt{tijojames@gmail.com}} 
\and 
Sandi Klav\v{z}ar$^{c,d,e,}$\thanks{Email: \texttt{sandi.klavzar@fmf.uni-lj.si}}
\and
Latha S Nair$^{a,}$\thanks{Email: \texttt{lathavichattu@gmail.com}} 
}
\maketitle

\begin{center}
$^a$ Department of Mathematics, Mar Athanasius College, Kothamangalam, India\\
	\medskip
$^b$ Department of Mathematics, Pavanatma College, Murickassery, India\\
\medskip
$^c$ Faculty of Mathematics and Physics, University of Ljubljana, Slovenia \\
\medskip

$^{d}$ Institute of Mathematics, Physics and Mechanics, Ljubljana, Slovenia \\
\medskip

$^e$ Faculty of Natural Sciences and Mathematics, University of Maribor, Slovenia \\
\medskip
\end{center}

\begin{abstract}
The Maker-Breaker total domination number, $\gamma_{\rm MBT}(G)$, of a graph $G$ is introduced as the minimum number of moves of Dominator to win the Maker-Breaker total domination game, provided that he has a winning strategy and is the first to play. The Staller-start Maker-Breaker total domination number, $\gamma_{\rm MBT}'(G)$, is defined analogously for the game in which Staller starts. Upper and lower bounds on $\gamma_{\rm MBT}(G)$ and on $\gamma_{\rm MBT}'(G)$ are provided and demonstrated to be sharp. It is proved that for any pair of integers $(k,\ell)$ with $2\leq k\leq \ell$, (i) there exists a connected graph $G$ with $\gamma_{\rm MB}(G)=k$ and $\gamma_{\rm MBT}(G)=\ell$, (ii) there exists a connected graph $G'$ with $\gamma_{\rm MB}'(G')=k$ and $\gamma_{\rm MBT}'(G')=\ell$, and (iii) there there exists a connected graph $G''$ with $\gamma_{\rm MBT}(G'')=k$ and $\gamma_{\rm MBT}'(G'')=\ell$. Here, $\gamma_{\rm MB}$ and $\gamma_{\rm MB}'$ are corresponding invariants for the Maker-Breaker domination game. 
\end{abstract}


\medskip\noindent
{\bf AMS Subj.\ Class.\ (2020)}: 05C57, 05C69

 \section{Introduction}
\label{sec: intro}

Maker-Breaker game is a two-player game played on an arbitrary hypergraph by Maker and Breaker. The game belongs to the larger family of Maker-Breaker positional games introduced by Hales and Jewett~\cite{HJ-1963}, and later by Erd\H{o}s and Selfridge~\cite{erdos-1973}. The game has been the subject of extensive study, both in general settings and specific cases, as detailed in the book~\cite{hefetz-2014}. See also the recent paper~\cite{rahman-2023} and references therein.

Duch\^{e}ne, Gledel, ~Parreau, and Renault~\cite{duchene-2020} introduced a variation of the Maker-Breaker game called the {\em Maker-Breaker domination game} (or {\em MBD game} for short). The {\em Maker-Breaker domination game} is played on a graph $G=(V(G), E(G))$ by two players named Dominator and Staller. These names align with the terminology used in the well-studied domination game~\cite{bresar-2010, book-2021}. The players alternately select unplayed vertices of $G$. Dominator's goal is to occupy all vertices of some dominating set of $G$, while Staller prevents it from happening. Recent research on the game includes~\cite{bagan-2025, divakaran-2025}. 

The total domination game from~\cite{henning-2015} has received almost as much attention as the domination game, see for instance~\cite{irsic-2019, jiang-2019, worawannotai-2024}. We would particularly like to highlight that recently Portier and Versteegen~\cite{portier-2025} proved the $3/4$-conjecture posed in~\cite{henning-2017}. Motivated by the total domination game, the {\em Maker-Breaker total domination game} ({\em MBTD game} for short) was introduced in~\cite{gledel-2020}. The game is played just like the Maker-Breaker domination game, except that Dominator aims to occupy all the vertices of a total dominating set of $G$. In~\cite{gledel-2020}, the outcome of the MBTD game was determined for cacti, and for Cartesian products of paths and cycles, while in~\cite{forcan-2022} the focus was on cubic graphs. Let us also mention here a related game domination subdivision number studied by Favaron, Karami, and Sheikholeslami~\cite{favaron-2015}. 

A {\rm D-game} is a MBD game or a MBTD in which Dominator takes the first turn, otherwise we speak of an {\em S-game}. There are four graph invariants associated with the MBD game~\cite{bujtas-2024, gledel-2019}. The {\em Maker-Breaker domination number}, $\gmb(G)$, represents the minimum number of moves required for Dominator to win the D-game on $G$, assuming an optimal play by both players. Similarly, the {\em Staller-Maker-Breaker domination number}, $\gsmb(G)$, represents the minimum number of moves required for Staller to win the D-game. The parameters $ \gmb'(G)$ and $\gsmb'(G)$ are defined for the S-game. The study in~\cite{bujtas-2024} focuses on $\gsmb(G)$ and $\gsmb'(G)$, providing, among other results, exact formulas for $\gsmb(G)$ when $G$ is a path. In~\cite{bujtas-2023}, trees $T$ with $\gsmb(T)=k$ were characterized for every positive integer $k$, and exact formulas for $\gsmb(G)$ and $\gsmb'(G)$ were derived for caterpillars. The main result of~\cite{forcan-2023} determined $\gmb$ and $\gmb'$ for Cartesian products of $K_2$ and a path. Additionally, in~\cite{dokyeesun-2024+}, the MBD game was further explored on Cartesian products of paths, stars, and complete bipartite graphs.

In this paper we introduce the Maker-Breaker total domination number $\gmbt(G)$ and the the Staller-start Maker-Breaker total domination number $\gmbt'(G)$ of a graph $G$. In Section~\ref{sec:bounds}, we provide upper and lower bounds on $\gmbt$ and on $\gmbt'$, and demonstrate their sharpness. Then, in Section~\ref{sec:realization}, we relate Maker-Breaker domination numbers with Maker-Breaker total domination numbers. We prove that for any pair of integers $(k,\ell)$ with $2\leq k\leq \ell$, (i) there exists a connected graph $G$ with $\gmb(G)=k$ and $\gmbt(G)=\ell$, (ii) there exists a connected graph $G'$ with $\gmb'(G')=k$ and $\gmbt'(G')=\ell$, and (iii) there there exists a connected graph $G''$ with $\gmbt(G'')=k$ and $\gmbt'(G'')=\ell$. In the rest of the section, additional definitions and concepts needed are given. 

For $k\in \mathbb{N}$ we will use the notation $[k] =\{1, \ldots, k\}$. Let $G=(V(G),E(G))$ be a graph. The order of $G$ is represented as $n(G)$. A {\em dominating set} of a graph $G$ is a subset $D \subseteq V(G)$ such that every vertex in $V(G)\setminus D$ has at least one neighbor in $D$, and $D$ is a {\em total dominating set} if every vertex in $V(G)$ has a neighbor in $D$. The {\em domination number} $\gamma(G)$ is the minimum cardinality of a dominating set in $G$, and the {\em total domination number} $\gamma_t(G)$ is the minimum cardinality of a total dominating set in $G$. A subset of pairs $\{(u_1, v_1), \ldots, (u_k, v_k)\}$ of vertices in $V(G)$ is a {\em pairing total dominating set} if all the vertices are distinct and for any selection of vertices $x_i\in \{u_i,v_i\}$, $i\in [k]$, the set $\{x_1, \dots, x_k\}$ is a total dominating set of $G$. 

The {\em outcome} $o(G)$ of the MBTD game played on $G$ is 
\begin{itemize}
\item $\mathcal{D}$, if Dominator has a winning strategy no matter who starts the game; 
\item $\mathcal{S}$, if Staller has a winning strategy no matter who starts the game; 
\item $\mathcal{N}$, if the first player has a winning strategy.
\end{itemize}
In an optimal strategy of Dominator (resp. Staller) in the MBTD game it is never an advantage
for him (resp.\ for her) to skip a move. This fact is known as {\bf No-Skip-Lemma}~\cite[Lemma 2.1]{gledel-2020}.

The {\em Maker-Breaker total domination number}, $\gmbt(G)$, of a graph $G$ is the minimum number of moves of Dominator to win the  MBTD game, provided that he has a winning strategy and is the first to play.  If Dominator has no winning strategy in the D-game, then set $\gmbt(G)= \infty$. The {\em Staller-start Maker-Breaker total domination number},  $\gsmb'(G)$,  is defined analogously for the S-game. 

\section{Bounds on the MBTD numbers}
\label{sec:bounds}

In this section we demonstrate that straightforward upper and lower bounds on $\gmbt$ and on $\gmbt'$ are sharp. The result dealing with upper bounds reads as follows. 

\begin{theorem} \label{thm:bound}
Let G be a graph.
\begin{itemize}
\item[(i)] If $\gmbt(G) <\infty$, then $\gmbt(G) \leq \left \lceil\frac{n(G)}{2}\right \rceil$.
\item[(ii)] If $\gmbt'(G) < \infty$, then $ \gmbt'(G) \leq \left \lfloor\frac{n(G)}{2}\right \rfloor$.
\end{itemize}
Moreover, for any $k > 1$, there is a connected graph $H$ with  $\gmbt(H) = k = \left\lceil\frac{n(H)}{2}\right \rceil$, and  there is a connected graph $H'$ with  $\gmbt'(H')=k=\left\lfloor\frac{n(H')}{2}\right \rfloor$.
 \end{theorem}
\proof
Let $G$ be a graph with $\gmbt(G)<\infty$. Since Dominator and Staller play alternately and Dominator starts the game, he can select at most $\left \lceil\frac{n(G)}{2}\right \rceil$ vertices, hence $\gmbt(G) \leq \left \lceil\frac{n(G)}{2}\right \rceil$. Analogous argument applies to the inequality (ii). 

To prove  sharpness of the bounds we give four constructions based on the parity of $k$, and whether the D-game or the S-game is played. 

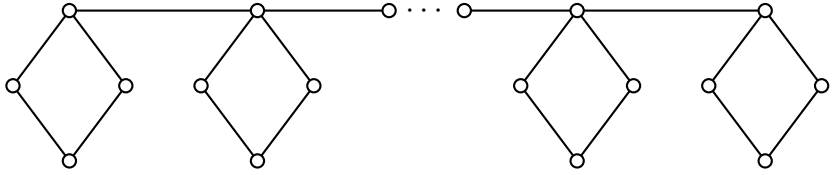
\begin{figure}[ht!]
\begin{center}
\begin{tikzpicture}[scale=0.5,style=thick,x=1cm,y=1cm]
\def\vr{5pt}
\begin{scope}[xshift=0cm, yshift=0cm] 
\coordinate(v_{11}) at (-2,12);
\coordinate(v_{12}) at (-.5,14);
\coordinate(v_{13}) at (-.5,10);
\coordinate(v_{14}) at (1,12);
\draw (v_{11}) -- (v_{12});
\draw (v_{11}) -- (v_{13});
\draw (v_{12}) -- (v_{14});
\draw (v_{13}) -- (v_{14});
\coordinate(v_{21}) at (3,12);
\coordinate(v_{22}) at (4.5,14);
\coordinate(v_{23}) at (4.5,10);
\coordinate(v_{24}) at (6,12);
\draw (v_{21}) -- (v_{22});
\draw (v_{21}) -- (v_{23});
\draw (v_{22}) -- (v_{24});
\draw (v_{23}) -- (v_{24});
\coordinate(x_1) at (8,14);
\draw (v_{12}) -- (v_{22});
\draw (v_{22}) -- (x_1);
\draw(v_{11})[fill=white] circle(\vr);
\draw(v_{12})[fill=white] circle(\vr);
\draw(v_{13})[fill=white] circle(\vr);
\draw(v_{14})[fill=white] circle(\vr);
\draw(v_{21})[fill=white] circle(\vr);
\draw(v_{22})[fill=white] circle(\vr);
\draw(v_{23})[fill=white] circle(\vr);
\draw(v_{24})[fill=white] circle(\vr);
\draw(x_1)[fill=white] circle(\vr);
\end{scope}
\begin{scope}[xshift=4cm, yshift=0cm] 
\node at (5,14) {$\cdots$};
\end{scope}
\begin{scope}[xshift=10cm, yshift=0cm] 
\coordinate(x_2) at (0,14);
\coordinate(v_{l-1 1}) at (1.5,12);
\coordinate(v_{l-1 2}) at (3,14);
\coordinate(v_{l-1 3}) at (3,10);
\coordinate(v_{l-1 4}) at (4.5,12);
\draw (v_{l-1 1}) -- (v_{l-1 2});
\draw (v_{l-1 1}) -- (v_{l-1 3});
\draw (v_{l-1 2}) -- (v_{l-1 4});
\draw (v_{l-1 3}) -- (v_{l-1 4});


\coordinate(v_{l 1}) at (6.5,12);
\coordinate(v_{l 2}) at (8,14);
\coordinate(v_{l 3}) at (8,10);
\coordinate(v_{l 4}) at (9.5,12);

\draw (v_{l 1}) -- (v_{l 2});
\draw (v_{l 1}) -- (v_{l 3});
\draw (v_{l 2}) -- (v_{l 4});
\draw (v_{l 3}) -- (v_{l 4});


\draw (v_{l-1 2}) -- (v_{l 2});
\draw (v_{l-1 2}) -- (x_2);
\draw(v_{l-1 1})[fill=white] circle(\vr);
\draw(v_{l-1 2})[fill=white] circle(\vr);
\draw(v_{l-1 3})[fill=white] circle(\vr);
\draw(v_{l-1 4})[fill=white] circle(\vr);
\draw(v_{l 1})[fill=white] circle(\vr);
\draw(v_{l 2})[fill=white] circle(\vr);
\draw(v_{l 3})[fill=white] circle(\vr);
\draw(v_{l 4})[fill=white] circle(\vr);
\draw(x_2)[fill=white] circle(\vr);

\end{scope}
\end{tikzpicture}
\caption{Graph $G_\ell$}
\label{fig: union of C4}
\end{center}
\end{figure}

\medskip\noindent
{\bf{Case 1}}: D-game, $k$ even.\\
Let  $k = 2\ell$ and let $G_\ell$ be the graph of order $4\ell$ obtained from $P_\ell$ by respectively attaching a copy of $C_4$ to each of its vertices, see Fig.~\ref{fig: union of C4}. 

Since the vertex sets of the $\ell$ copies of $C_4$ partition $V(G_\ell)$ and since $o(C_4) = \mathcal{D}$, we can apply~\cite[Corollary 2.2(ii)]{gledel-2020} to conclude that  $o(G_\ell) = \mathcal{D}$. It remains to verify that $\gmbt(G_\ell)=\left \lceil\frac{n(G_\ell)}{2}\right \rceil = 2\ell = k$. By the upper bound (i) we have  $\gmbt(G_\ell)\leq \left \lceil\frac{n(G)}{2}\right \rceil = 2\ell = k$.  On the other hand, every total dominating set of $G_\ell$ necessarily contains at least two vertices from each of the $\ell$ $4$-cycles, hence $\gamma_t(G_\ell) \ge  2\ell = k$. So $\gamma_t(G_\ell) = 2\ell = k$, hence $\gmbt(G_\ell) \ge k$, and we can conclude that $\gmbt(G_\ell) = k$.
 
\medskip\noindent
{\bf{Case 2}}: D-game, $k$ odd.\\ 
Let $k=2\ell+1$, and let $G_\ell'$ be the graph constructed just as $G_\ell$, except that at the first vertex of the base path the graph as shown in Fig.~\ref{fig:union of $C_4$ and Twin $C_5$} is attached (and not $C_4$ as in the construction of $G_\ell$). 

\begin{figure}[ht!]
\begin{center}
\begin{tikzpicture}[scale=0.5,style=thick,x=1cm,y=1cm]
\def\vr{5pt}
\begin{scope}[xshift=0cm, yshift=0cm] 
\coordinate(v_{11}) at (-2,12);
\coordinate(v_{12}) at (-.5,14);
\coordinate(v_{13}) at (-.5,10);
\coordinate(v_{14}) at (1,12);
\coordinate(v_{15}) at (-2,8);
\coordinate(v_{16}) at (1,8);
\draw (v_{11}) -- (v_{12});
\draw (v_{11}) -- (v_{13});
\draw (v_{12}) -- (v_{14});
\draw (v_{13}) -- (v_{14});
\draw (v_{11}) -- (v_{15});
\draw (v_{14}) -- (v_{16});
\draw (v_{15}) -- (v_{16});
\coordinate(v_{21}) at (3,12);
\coordinate(v_{22}) at (4.5,14);
\coordinate(v_{23}) at (4.5,10);
\coordinate(v_{24}) at (6,12);
\draw (v_{21}) -- (v_{22});
\draw (v_{21}) -- (v_{23});
\draw (v_{22}) -- (v_{24});
\draw (v_{23}) -- (v_{24});
\node at (0.3,12.0) {$x$};
\coordinate(x_1) at (8,14);
\draw (v_{12}) -- (v_{22});
\draw (v_{22}) -- (x_1);
\draw(v_{21})[fill=white] circle(\vr);
\draw(v_{22})[fill=white] circle(\vr);
\draw(v_{23})[fill=white] circle(\vr);
\draw(v_{24})[fill=white] circle(\vr);
\draw(v_{11})[fill=white] circle(\vr);
\draw(v_{12})[fill=white] circle(\vr);
\draw(v_{13})[fill=white] circle(\vr);
\draw(v_{14})[fill=white] circle(\vr);
\draw(v_{15})[fill=white] circle(\vr);
\draw(v_{16})[fill=white] circle(\vr);
\draw(x_1)[fill=white] circle(\vr);
\end{scope}
\begin{scope}[xshift=4cm, yshift=0cm] 
\node at (5,14) {$\cdots$};
\end{scope}
\begin{scope}[xshift=10cm, yshift=0cm] 
\coordinate(x_2) at (0,14);
\coordinate(v_{l-1 1}) at (1.5,12);
\coordinate(v_{l-1 2}) at (3,14);
\coordinate(v_{l-1 3}) at (3,10);
\coordinate(v_{l-1 4}) at (4.5,12);
\draw (v_{l-1 1}) -- (v_{l-1 2});
\draw (v_{l-1 1}) -- (v_{l-1 3});
\draw (v_{l-1 2}) -- (v_{l-1 4});
\draw (v_{l-1 3}) -- (v_{l-1 4});
\coordinate(v_{l 1}) at (6.5,12);
\coordinate(v_{l 2}) at (8,14);
\coordinate(v_{l 3}) at (8,10);
\coordinate(v_{l 4}) at (9.5,12);
\draw (v_{l 1}) -- (v_{l 2});
\draw (v_{l 1}) -- (v_{l 3});
\draw (v_{l 2}) -- (v_{l 4});
\draw (v_{l 3}) -- (v_{l 4});
\draw (v_{l-1 2}) -- (v_{l 2});
\draw (v_{l-1 2}) -- (x_2);
\draw(v_{l-1 1})[fill=white] circle(\vr);
\draw(v_{l-1 2})[fill=white] circle(\vr);
\draw(v_{l-1 3})[fill=white] circle(\vr);
\draw(v_{l-1 4})[fill=white] circle(\vr);
\draw(v_{l 1})[fill=white] circle(\vr);
\draw(v_{l 2})[fill=white] circle(\vr);
\draw(v_{l 3})[fill=white] circle(\vr);
\draw(v_{l 4})[fill=white] circle(\vr);
\draw(x_2)[fill=white] circle(\vr);
\end{scope}
\end{tikzpicture}
\caption{Graph $G_\ell'$}
\label{fig:union of $C_4$ and Twin $C_5$}
\end{center}
\end{figure}

Let Dominator start the game by selecting the vertex $x$, see Fig.~\ref{fig:union of $C_4$ and Twin $C_5$}. Then a simple analysis reveals that Dominator has a winning strategy. By (i), 
$$\gmbt(G_\ell') \leq \left \lceil\frac{n(G_\ell')}{2}\right \rceil= \left \lceil\frac{4\ell + 2}{2}\right \rceil = 2\ell + 1 = k.$$
On the other hand, $\gamma_t(G_\ell') = 3 + 2(\ell-1) = 2\ell+1 = k$. Therefore, $\gmbt(G_\ell') = k$.

\medskip\noindent
{\bf{Case 3}}: S-game, $k$ even.\\
In this case, the graph $G_\ell$ from Case~1, where $k = 2\ell$, gives the required conclusion. The arguments are parallel, in particular, the bound (ii) is the same as the bound (i) because $k$ is even. 
 
\medskip\noindent
{\bf{Case 4}}: S-game, $k$ odd.\\
Let $k=2\ell+1$, and let $G_\ell''$ be the graph constructed just as the graphs $G_\ell$ and $G_\ell'$, except that at the first vertex of the base path the graph $H$ as shown in Fig.~\ref{fig:union of C4 and X89} is attached.

\begin{figure}[ht!]
\begin{center}
\begin{tikzpicture}[scale=0.5,style=thick,x=1cm,y=1cm]
\def\vr{5pt}
\begin{scope}[xshift=0cm, yshift=0cm] 
\coordinate(v_{1}) at (-2,12);
\coordinate(v_{2}) at (-.5,14);
\coordinate(v_{3}) at (-.5,10);
\coordinate(v_{4}) at (-.5,12);
\coordinate(v_{5}) at (1,12);
\coordinate(v_{6}) at (-2,8);
\coordinate(v_{7}) at (1,8);

\draw (v_{1}) -- (v_{2});
\draw (v_{1}) -- (v_{3});
\draw (v_{1}) -- (v_{4});
\draw (v_{2}) -- (v_{5});
\draw (v_{3}) -- (v_{5});
\draw (v_{4}) -- (v_{5});
\draw (v_{1}) -- (v_{6});
\draw (v_{5}) -- (v_{7});
\draw (v_{6}) -- (v_{7});
\draw (v_{3}) -- (v_{7});
\draw (v_{3}) -- (v_{6});

\node at (-2.7,12.2) {$x_{2}$};
\node at (-.5,14.6) {$x_{5}$};
\node at (-.5,10.8)  {$x_{3}$};
\node at (1.7,12.2) {$x_{1}$};
\node at (-2.7,8.1) {$x_{4}$};
\node at (-0.7,12.6) {$x_{6}$};

\coordinate(v_{21}) at (3,12);
\coordinate(v_{22}) at (4.5,14);
\coordinate(v_{23}) at (4.5,10);
\coordinate(v_{24}) at (6,12);

\draw (v_{21}) -- (v_{22});
\draw (v_{21}) -- (v_{23});
\draw (v_{22}) -- (v_{24});
\draw (v_{23}) -- (v_{24});


\coordinate(x_1) at (8,14);
\draw (v_{12}) -- (v_{22});
\draw (v_{22}) -- (x_1);
\draw(v_{21})[fill=white] circle(\vr);
\draw(v_{22})[fill=white] circle(\vr);
\draw(v_{23})[fill=white] circle(\vr);
\draw(v_{24})[fill=white] circle(\vr);
\draw(v_{1})[fill=white] circle(\vr);
\draw(v_{2})[fill=white] circle(\vr);
\draw(v_{3})[fill=white] circle(\vr);
\draw(v_{4})[fill=white] circle(\vr);
\draw(v_{5})[fill=white] circle(\vr);
\draw(v_{6})[fill=white] circle(\vr);
\draw(v_{7})[fill=white] circle(\vr);
\draw(x_1)[fill=white] circle(\vr);

\end{scope}

\begin{scope}[xshift=4cm, yshift=0cm] 

\node at (5,14) {$\cdots$};

\end{scope}

\begin{scope}[xshift=10cm, yshift=0cm] 

\coordinate(x_2) at (0,14);

\coordinate(v_{l-1 1}) at (1.5,12);
\coordinate(v_{l-1 2}) at (3,14);
\coordinate(v_{l-1 3}) at (3,10);
\coordinate(v_{l-1 4}) at (4.5,12);

\draw (v_{l-1 1}) -- (v_{l-1 2});
\draw (v_{l-1 1}) -- (v_{l-1 3});
\draw (v_{l-1 2}) -- (v_{l-1 4});
\draw (v_{l-1 3}) -- (v_{l-1 4});


\coordinate(v_{l 1}) at (6.5,12);
\coordinate(v_{l 2}) at (8,14);
\coordinate(v_{l 3}) at (8,10);
\coordinate(v_{l 4}) at (9.5,12);

\draw (v_{l 1}) -- (v_{l 2});
\draw (v_{l 1}) -- (v_{l 3});
\draw (v_{l 2}) -- (v_{l 4});
\draw (v_{l 3}) -- (v_{l 4});


\draw (v_{l-1 2}) -- (v_{l 2});
\draw (v_{l-1 2}) -- (x_2);
\draw(v_{l-1 1})[fill=white] circle(\vr);
\draw(v_{l-1 2})[fill=white] circle(\vr);
\draw(v_{l-1 3})[fill=white] circle(\vr);
\draw(v_{l-1 4})[fill=white] circle(\vr);
\draw(v_{l 1})[fill=white] circle(\vr);
\draw(v_{l 2})[fill=white] circle(\vr);
\draw(v_{l 3})[fill=white] circle(\vr);
\draw(v_{l 4})[fill=white] circle(\vr);
\draw(x_2)[fill=white] circle(\vr);

\end{scope}

\end{tikzpicture}
\caption{Graph $G_\ell''$}
\label{fig:union of C4 and X89}
\end{center}
\end{figure}

To see that Dominator wins the S-game, note that as before, he has an easy control in the game in the $4$-cycles. As for the subgraph $H$, assume that at some point of the game, possible the first move of it, Staller has selected the vertex $x_1$. Then Dominator is forced to play $x_2$. After that, Staller selects $x_3$ which in turn forces Dominator to select $x_4$. Now the only vertex in $H$ which is not yet totally dominated is $x_1$. By selecting $x_5$ or $x_6$ in the next move, Dominator totally dominates $H$ in three moves. We can argue similarly that Dominator can win if Staller selects any other vertex of $H$, notably $x_3$. Moreover, the above strategy of Staller implies that $\gmbt'(G_\ell'') \geq 3 + 2(\ell - 1)$, so that $\gmbt'(G_\ell'') \ge 2\ell+1 = k$. On the other hand, (ii) yields $$\gmbt'(G_\ell'') \leq \left \lfloor\frac{n(G_\ell'')}{2}\right \rfloor= \left \lfloor\frac{4\ell+3}{2}\right \rfloor=2\ell+1 = k\,,$$
and we are done. 
\qed

We now turn our attention to lower bounds. Let $G$ be a graph. A winning set for Dominator in a Maker-Breaker total domination game is a total dominating set, hence, $\gamma_t(G)\leq\gmbt(G)$. Since a D-game can be viewed as an S-game in which Staller has skipped her first move, the No-Skip Lemma implies $\gmbt(G)\leq\gmbt'(G)$. Moreover, since every total dominating set is a dominating set, we also infer that $\gmb(G)\leq\gmbt(G)$. These observations can be summarized as follows: 
 
\begin{equation}
\label{eq:basic-inequalities}
\max \{\gamma_t(G), \gmb(G)\} \leq \gmbt(G) \leq \gmbt'(G)\,.
\end{equation}

In view of \eqref{eq:basic-inequalities}, we now prove the following result which demonstrate that simultaneous equalities can hold in these inequalities. 

\begin{theorem}\label{equal graphs for D game}
If $k\ge 2$, then the following hold. 
\begin{enumerate}
\item[(i)] There exists a connected graph $G$ with $\gamma(G) = \gamma_t(G) = \gmb(G)=\gmbt(G)$.
\item[(ii)] There exists a connected graph $H$ with $\gamma(H)=\gamma_t(H)=\gmb'(H)=\gmbt'(H)$.
\end{enumerate}
\end{theorem}

\proof
(i) Let $G_{k,n}$, $k\ge 2$, $n\ge 4$, be the graph constructed in the following way. Take the disjoint union of $K_n$ and $k-1$ copies of $K_n-e$. Then add one more vertex $u$, and add edges between $u$ and all the vertices of $K_n$, and between $u$ and each pair of non-adjacent vertices of each $K_n-e$. See Fig.~\ref{fig:equal graph}.

\begin{figure}[ht!]
\begin{center}
\begin{tikzpicture}[scale=0.45,style=thick,x=1cm,y=1cm]
\def\vr{6pt}
\begin{scope}[xshift=0cm, yshift=0cm] 
\coordinate(y_1) at (7,13);
\coordinate(y_2) at (13,13);
\coordinate(y_3) at (7,17);
\coordinate(y_4) at (13,17);
\coordinate(u) at (10,9);

\node at (9.1,9.2) {$u$};

\draw (y_1) -- (y_2);
\draw (y_1) -- (y_3);
\draw (y_1) -- (y_4);
\draw (y_2) -- (y_3);
\draw (y_2) -- (y_4);
\draw (y_3) -- (y_4);

\draw (y_1) -- (u);
\draw (y_2) -- (u);
\draw (y_3) -- (u);
\draw (y_4) -- (u);

\draw(y_1)[fill=white] circle(\vr);
\draw(y_2)[fill=white] circle(\vr);
\draw(y_3)[fill=white] circle(\vr);
\draw(y_4)[fill=white] circle(\vr);

\coordinate(x_{11}) at (-3,1);
\coordinate(x_{12}) at (-3,4);
\coordinate(x_{13}) at (1,1);
\coordinate(x_{14}) at (1,4);

\draw (x_{12}) -- (u);
\draw (x_{14}) -- (u);
\draw (x_{11}) -- (x_{12});
\draw (x_{11}) -- (x_{13});
\draw (x_{11}) -- (x_{14});
\draw (x_{12}) -- (x_{13});
\draw (x_{13}) -- (x_{14});

\draw(x_{11})[fill=white] circle(\vr);
\draw(x_{12})[fill=white] circle(\vr);
\draw(x_{13})[fill=white] circle(\vr);
\draw(x_{14})[fill=white] circle(\vr);

\coordinate(x_{21}) at (4,1);
\coordinate(x_{22}) at (4,4);
\coordinate(x_{23}) at (8,1);
\coordinate(x_{24}) at (8,4);

\draw (x_{22}) -- (u);
\draw (x_{24}) -- (u);
\draw (x_{21}) -- (x_{22});
\draw (x_{21}) -- (x_{23});
\draw (x_{21}) -- (x_{24});
\draw (x_{22}) -- (x_{23});
\draw (x_{23}) -- (x_{24});

\draw(x_{21})[fill=white] circle(\vr);
\draw(x_{22})[fill=white] circle(\vr);
\draw(x_{23})[fill=white] circle(\vr);
\draw(x_{24})[fill=white] circle(\vr);


\node at (10,2.5) {$\cdots$};

\coordinate(x_{k-1 1}) at (12,1);
\coordinate(x_{k-1 2}) at (12,4);
\coordinate(x_{k-1 3}) at (16,1);
\coordinate(x_{k-1 4}) at (16,4);

\draw (x_{k-1 2}) -- (u);
\draw (x_{k-1 4}) -- (u);
\draw (x_{k-1 1}) -- (x_{k-1 2});
\draw (x_{k-1 1}) -- (x_{k-1 3});
\draw (x_{k-1 1}) -- (x_{k-1 4});
\draw (x_{k-1 2}) -- (x_{k-1 3});
\draw (x_{k-1 3}) -- (x_{k-1 4});

\draw(x_{k-1 1})[fill=white] circle(\vr);
\draw(x_{k-1 2})[fill=white] circle(\vr);
\draw(x_{k-1 3})[fill=white] circle(\vr);
\draw(x_{k-1 4})[fill=white] circle(\vr);


\coordinate(x_{k 1}) at (19,1);
\coordinate(x_{k 2}) at (19,4);
\coordinate(x_{k 3}) at (23,1);
\coordinate(x_{k 4}) at (23,4);

\draw (x_{k 2}) -- (u);
\draw (x_{k 4}) -- (u);
\draw (x_{k 1}) -- (x_{k 2});
\draw (x_{k 1}) -- (x_{k 3});
\draw (x_{k 1}) -- (x_{k 4});
\draw (x_{k 2}) -- (x_{k 3});
\draw (x_{k 3}) -- (x_{k 4});

\draw(x_{k 1})[fill=white] circle(\vr);
\draw(x_{k 2})[fill=white] circle(\vr);
\draw(x_{k 3})[fill=white] circle(\vr);
\draw(x_{k 4})[fill=white] circle(\vr);
\draw(u) [fill=white] circle(\vr);

\end{scope}
\end{tikzpicture}
    \caption{Graph $G_{k,4}$ with $\gamma(G_{k,4})=\gamma_t(G_{k,4})=\gmb(G_{k,4})=\gmbt(G_{k,4}) = k$}
    \label{fig:equal graph}
    \end{center}
\end{figure}

It is straightforward to see that $\gamma(G_{k,n}) = \gamma_t(G_{k,n})=k$. In particular, the vertex $u$ together with exactly one vertex of each of $K_n-e$ adjacent to $u$ form a $\gamma_t$-set. Consider now the D-game. Then  Dominator selects $u$ as his first move. After that his strategy is to select one vertex of each of $K_n-e$ adjacent to $u$. In this way Dominator wins in $k$ moves, and as $\gmbt(G_{k,n}) \ge \gamma_t(G_{k,n})$ we can conclude that $\gmbt(G_{k,n}) = k$. Consequently,  $\gmb(G_{k,n}) = k$ holds as well. 

(ii) Let $H_{k,n}$, $k\ge 2$, $n\ge 4$, be the graph constructed similarly as the graph $G_{k,n}$ from Fig.~\ref{fig:equal graph}. The only difference is that instead of one additional vertex $u$ we take two additional vertices $u$ and $v$ which are then adjacent to all the vertices to which $u$ is adjacent in $G_{k,n}$, see Fig.~\ref{fig:equal graph for S game}. 

\begin{figure}[ht!]
\begin{center}
\begin{tikzpicture}[scale=0.45,style=thick,x=1cm,y=1cm]
\def\vr{6pt}
\begin{scope}[xshift=0cm, yshift=0cm] 
\coordinate(y_1) at (7,15);
\coordinate(y_2) at (13,15);
\coordinate(y_3) at (7,19);
\coordinate(y_4) at (13,19);
\coordinate(u) at (9,11);
\coordinate (v) at (11,11);

\node at (8,11) {$u$};
\node at (12,11) {$v$};

\draw (y_1) -- (y_2);
\draw (y_1) -- (y_3);
\draw (y_1) -- (y_4);
\draw (y_2) -- (y_3);
\draw (y_2) -- (y_4);
\draw (y_3) -- (y_4);

\draw [thick, blue] (y_1) -- (u);
\draw [thick, blue] (y_2) -- (u);
\draw [thick, blue](y_3) -- (u);
\draw [thick, blue](y_4) -- (u);
\draw [thick, red](y_1) -- (v);
\draw [thick, red](y_2) -- (v);
\draw [thick, red](y_3) -- (v);
\draw [thick, red](y_4) -- (v);

\draw(y_1)[fill=white] circle(\vr);
\draw(y_2)[fill=white] circle(\vr);
\draw(y_3)[fill=white] circle(\vr);
\draw(y_4)[fill=white] circle(\vr);

\coordinate(x_{11}) at (-3,1);
\coordinate(x_{12}) at (-3,4);
\coordinate(x_{13}) at (1,1);
\coordinate(x_{14}) at (1,4);

\draw [thick, blue](x_{12}) -- (u);
\draw [thick, blue](x_{14}) -- (u);
\draw [thick, red](x_{12}) -- (v);
\draw [thick, red](x_{14}) -- (v);
\draw (x_{11}) -- (x_{12});
\draw (x_{11}) -- (x_{13});
\draw (x_{11}) -- (x_{14});
\draw (x_{12}) -- (x_{13});
\draw (x_{13}) -- (x_{14});

\draw(x_{11})[fill=white] circle(\vr);
\draw(x_{12})[fill=white] circle(\vr);
\draw(x_{13})[fill=white] circle(\vr);
\draw(x_{14})[fill=white] circle(\vr);


\coordinate(x_{21}) at (4,1);
\coordinate(x_{22}) at (4,4);
\coordinate(x_{23}) at (8,1);
\coordinate(x_{24}) at (8,4);

\draw [thick, blue](x_{22}) -- (u);
\draw [thick, blue](x_{24}) -- (u);
\draw [thick, red](x_{22}) -- (v);
\draw [thick, red](x_{24}) -- (v);
\draw (x_{21}) -- (x_{22});
\draw (x_{21}) -- (x_{23});
\draw (x_{21}) -- (x_{24});
\draw (x_{22}) -- (x_{23});
\draw (x_{23}) -- (x_{24});

\draw(x_{21})[fill=white] circle(\vr);
\draw(x_{22})[fill=white] circle(\vr);
\draw(x_{23})[fill=white] circle(\vr);
\draw(x_{24})[fill=white] circle(\vr);

\node at (10,2.5) {$\cdots$};

\coordinate(x_{k-1 1}) at (12,1);
\coordinate(x_{k-1 2}) at (12,4);
\coordinate(x_{k-1 3}) at (16,1);
\coordinate(x_{k-1 4}) at (16,4);

\draw(x_{k-1 1})[fill=white] circle(\vr);
\draw(x_{k-1 2})[fill=white] circle(\vr);
\draw(x_{k-1 3})[fill=white] circle(\vr);
\draw(x_{k-1 4})[fill=white] circle(\vr);

\draw [thick, blue](x_{k-1 2}) -- (u);
\draw [thick, blue](x_{k-1 4}) -- (u);
\draw [thick, red](x_{k-1 2}) -- (v);
\draw [thick, red](x_{k-1 4}) -- (v);
\draw (x_{k-1 1}) -- (x_{k-1 2});
\draw (x_{k-1 1}) -- (x_{k-1 3});
\draw (x_{k-1 1}) -- (x_{k-1 4});
\draw (x_{k-1 2}) -- (x_{k-1 3});
\draw (x_{k-1 3}) -- (x_{k-1 4});

\draw(x_{k-1 1})[fill=white] circle(\vr);
\draw(x_{k-1 2})[fill=white] circle(\vr);
\draw(x_{k-1 3})[fill=white] circle(\vr);
\draw(x_{k-1 4})[fill=white] circle(\vr);


\coordinate(x_{k 1}) at (19,1);
\coordinate(x_{k 2}) at (19,4);
\coordinate(x_{k 3}) at (23,1);
\coordinate(x_{k 4}) at (23,4);

\draw [thick, blue](x_{k 2}) -- (u);
\draw [thick, blue](x_{k 4}) -- (u);
\draw [thick, red](x_{k 2}) -- (v);
\draw [thick, red](x_{k 4}) -- (v);
\draw (x_{k 1}) -- (x_{k 2});
\draw (x_{k 1}) -- (x_{k 3});
\draw (x_{k 1}) -- (x_{k 4});
\draw (x_{k 2}) -- (x_{k 3});
\draw (x_{k 3}) -- (x_{k 4});

\draw(x_{k 1})[fill=white] circle(\vr);
\draw(x_{k 2})[fill=white] circle(\vr);
\draw(x_{k 3})[fill=white] circle(\vr);
\draw(x_{k 4})[fill=white] circle(\vr);
\draw(u) [fill=white] circle(\vr);
\draw(v) [fill=white] circle(\vr);
\end{scope}

\end{tikzpicture}
     \caption{Graph $H_{k,n}$ with $\gamma(H_{k,n}) = \gamma_t(H) = \gmb'(H) = \gmbt'(H) = k$}
    \label{fig:equal graph for S game}
    \end{center}
\end{figure}

Let $u_i, v_i$, $i\in  [k-1]$, be the two non-adjacent vertices of the $i^{\rm th}$ copy of $K_n-e$. Then $\{(u,v), (u_1, v_1), \ldots, (u_{k-1}, v_{k-1})\}$ is a total pairing dominating set. Therefore Dominator can win the game in $k$ moves and since $\gamma_t(H_{k,n}) = k$, we get $\gmbt'(H)=k$. Hence the result. 
\qed

\section{Realization of parameters}
\label{sec:realization}

In this section we relate Maker-Breaker domination numbers with Maker-Breaker total domination numbers and demonstrate that in all non-trivial cases the difference can be arbitrarily large.

If $\gmb(G)=1$, then $\gamma(G)=1$. Hence, if $n(G)\ge 3$, then $\gmbt(G) = 2$. Similarly, if $\gmb'(G)=1$, then $G$ contains at least two vertices of degree $n(G)-1$. Therefore, if $n(G)\ge 4$, then $\gmbt'(G) = 2$. In this section, we build on these observations and show that for any pair of integers $(k,\ell)$ with $2\leq k\leq \ell$, there exists a connected graph $G$ with $\gmb(G)=k$ and $\gmbt(G)=\ell$, as well as a connected graph $H$ with $\gmb'(G)=k$ and $\gmbt'(G)=\ell$. Several different constructions are needed for this task because the case $k=2$ behaves special.  

\begin{theorem}
\label{thm:MB-MBT}
  For any two integers $k$ and $\ell$ with $2\leq k\leq  \ell$,  there exists a connected graph $G$ with $\gmb(G)=k$ and $\gmbt(G)=\ell$. 
\end{theorem}
\proof
Clearly, $\gamma(C_4)=\gmb(C_4)=2$  and $\gmbt(C_4)=2$. Assume in the rest that $\ell\geq 3$. We consider two cases. 

\begin{figure}[ht!]
\begin{center}
\begin{tikzpicture}[scale=0.45,style=thick,x=1cm,y=1cm]
\def\vr{6pt}
\begin{scope}[xshift=0cm, yshift=0cm] 

\coordinate(u) at (6,9);
\coordinate (v) at (14,9);
\node at (5,9) {$u$};
\node at (15,9) {$v$};
\draw (u) -- (v);

\coordinate(v_{21}) at (-3,16);
\coordinate(v_{22}) at (1,16);
\coordinate(v_{23}) at (1,19);
\node at (-4,16) {$v_{2,1}$};
\node at (0,16.4) {$v_{2,2}$};
\node at (1,19.8) {$v_{2,3}$};

\draw [thick, blue](v_{21}) -- (u);
\draw [thick, blue](v_{23}) -- (u);
\draw [thick, blue](v_{22}) -- (u);
\draw [thick, red](v_{21}) -- (v);
\draw [thick, red](v_{22}) -- (v);

\draw (v_{21}) -- (v_{22});
\draw (v_{21}) -- (v_{23});
\draw (v_{22}) -- (v_{23});

\coordinate(v_{31}) at (8,17);
\coordinate(v_{32}) at (8,20);
\coordinate(v_{33}) at (12,17);

\draw [thick, blue](v_{31}) -- (u);
\draw [thick, blue](v_{33}) -- (u);
\draw [thick, blue](v_{32}) -- (u);
\draw [thick, red](v_{31}) -- (v);
\draw [thick, red](v_{33}) -- (v);

\draw (v_{31}) -- (v_{32});
\draw (v_{31}) -- (v_{33});
\draw (v_{32}) -- (v_{33});

\node at (14.5,18) {$\cdots$};

\coordinate(v_{l-1 1}) at (19,16);
\coordinate(v_{l-1 2}) at (19,19);
\coordinate(v_{l-1 3}) at (23,16);

\node at (17.9,16.5) {$v_{\ell-1, 1}$};
\node at (24.5,16) {$v_{\ell-1, 2}$};
\node at (19,19.5) {$v_{\ell-1, 3}$};

\draw [thick, blue](v_{l-1 1}) -- (u);
\draw [thick, blue](v_{l-1 3}) -- (u);
\draw [thick, blue](v_{l-1 2}) -- (u);
\draw [thick, red](v_{l-1 1}) -- (v);
\draw [thick, red](v_{l-1 3}) -- (v);
\draw (v_{l-1 1}) -- (v_{l-1 2});
\draw (v_{l-1 1}) -- (v_{l-1 3});
\draw (v_{l-1 2}) -- (v_{l-1 3});

\coordinate(v_{1 1}) at (10.5,4);
\coordinate(v_{1 2}) at (7.5,1);
\coordinate(v_{1 3}) at (14,1);
\node at (11.7,4) {$v_{1, 1}$};
\node at (6,1) {$v_{1, 2}$};
\node at (15,1) {$v_{1, 3}$};

\draw (v_{1 1}) -- (v_{1 2});
\draw (v_{1 1}) -- (v_{1 3});
\draw (v_{1 2}) -- (v_{1 3});

\draw [thick, blue](v_{1 1}) -- (u);
\draw [thick, red](v_{1 1}) -- (v);

\draw(v_{1 1})[fill=white] circle(\vr);
\draw(v_{1 2})[fill=white] circle(\vr);
\draw(v_{1 3})[fill=white] circle(\vr);

\draw(v_{21})[fill=white] circle(\vr);
\draw(v_{22})[fill=white] circle(\vr);
\draw(v_{23})[fill=white] circle(\vr);

\draw(v_{31})[fill=white] circle(\vr);
\draw(v_{32})[fill=white] circle(\vr);
\draw(v_{33})[fill=white] circle(\vr);

\draw(v_{l-1 1})[fill=white] circle(\vr);
\draw(v_{l-1 2})[fill=white] circle(\vr);
\draw(v_{l-1 3})[fill=white] circle(\vr);
\draw(u)[fill=white] circle(\vr);
\draw(v)[fill=white] circle(\vr);
\end{scope}
\end{tikzpicture}
\caption{Graph $G_{2,\ell}$ with $\gmb(G_{2,\ell})=2$ and $\gmbt(G_{2,\ell}) = \ell$}
    \label{fig:real graph 1}
    \end{center}
\end{figure}
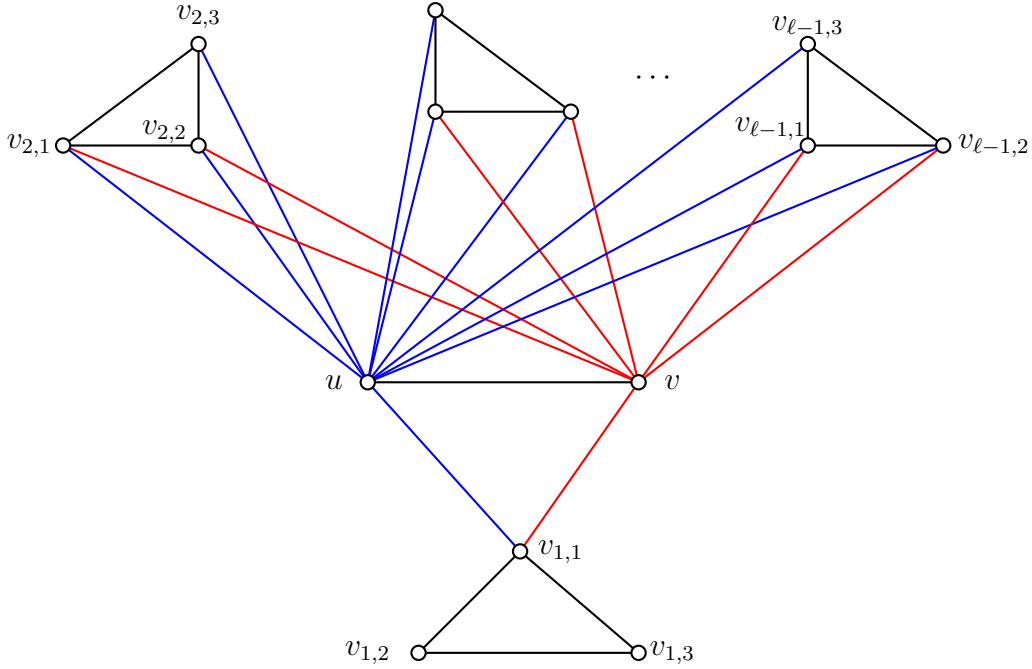

\medskip\noindent 
{\bf Case 1}: $k=2$. \\
Let $G_{2,\ell}$, $\ell \ge 3$, be the graph constructed from the disjoint union of $\ell-1$ triangles as follows. Let $v_{i,1}, v_{i,2}, v_{i,3}$ be the vertices of the $i^{\rm th}$ triangle, $i\in [\ell -1]$. Take two additional vertices $u$ and $v$ both adjacent to $v_{1,1}$. In addition, $u$ is adjacent to $v_{i,j}$  for $i \in \{2,3,\ldots, \ell-1\}$ and $j\in [3]$. Also, $v$ is adjacent to the vertices $v_{i,1}, v_{i,2}$ for $i\in \{2,3,\ldots, \ell-1\}$. See Fig.~\ref{fig:real graph 1}.

Consider a Maker-Breaker domination game on $G_{2, \ell}$. The sets $\{u, v_{1,1}\}$ and $\{u, v_{1,2}\}$ are dominating sets of $G_{2, \ell}$. Dominator first selects $u$, and then he can select either  $v_{1,1}$ or $v_{1,2}$ as his second move. Therefore, Dominator can finish the game in two moves and hence $\gmb(G_{2, \ell})=2$ because $G_{2, \ell}$ has no vertex adjacent to all the other vertices. 

Now we prove $\gmbt(G_{2, \ell}) =\ell$. To show that $\gmbt(G_{2, \ell})\leq \ell$, we need to describe a strategy for Dominator in which he can finish the game in at most $\ell$ moves. Dominator starts by selecting $v_{1,1}$. Then Staller replies by selecting $u$, for otherwise  Dominator can finish the game in two moves by playing $u$ as his second move. Then Dominator selects $v$ as his second move. The only remaining undominated vertices are $v_{2,3},v_{3,3},\dots, v_{\ell-1,3}$. So in the continuation of the game, Dominator selects either $v_{i,1}$ or $v_{i,2}$  for $i\in\{2, 3,\ldots, \ell-1\}$ which he can clearly achieve. Hence Dominator can finish the game in at most $2+\ell-2=\ell$ moves and so $\gmbt(G)\leq \ell$. 

To prove that $\gmbt(G_{2, \ell})\geq \ell$, we describe a strategy of Staller. Suppose that Dominator selects a vertex other than $v_{1,j}$, $j\in[3]$. In this case, Staller selects $v_{1,1}$ as her response and wins the game by playing as her next move $v_{1,2}$ or $v_{1,3}$. Hence Dominator's first move is a vertex of the first triangle. If he selects $v_{1,2}$ or $v_{1,3}$, say $v_{1,2}$, then Staller selects $v_{1,1}$. In this case, Dominator is forced next to select $v_{1,3}$. Now Staller selects $u$ as her second move which forces Dominator to play $v$, and at least $\ell-2$ more moves. Finally, if Dominator first plays $v_{1,1}$, then Staller replies by selecting $u$ as her first move, which forces Dominator to play at  least $\ell-1$ more moves to finish the game (where $v$ is his second move). Hence $\gmbt(G_{2, \ell})\geq \ell$ and we get $\gmbt(G_{2, \ell})= \ell$.

\medskip\noindent
{\bf Case 2}: $k\ge 3$. \\
Construct the graph $G_{k,\ell}$, $k\ge 3$, $\ell \ge 3$, as follows. Take $k+\ell-2$ triangles, where the $i^{\rm th}$ triangle, $i\in[k+\ell-2]$, has vertices $v_{i,j}$, $j \in [3]$. Then add two more vertices, $u$ and $v$. The vertex $u$ is adjacent to $v_{i,1}$, $i\in [k-1]$, and to $v_{i,j}$, $i\in \{k, k+1,\ldots, k+\ell-2\}$, $j\in [3]$. The vertex $v$ is adjacent to $v_{i,j}$, $i\in [k-1]$, $j\in [3]$, and to $v_{i,j}$, $i\in \{k, k+1,\ldots, k+\ell-2\}$, $j\in [2]$. See Fig.~\ref{fig:real graph 2}.

\begin{figure}[ht!]
\begin{center}
\begin{tikzpicture}[scale=0.45,style=thick,x=1cm,y=1cm]
\def\vr{6pt}
\begin{scope}[xshift=0cm, yshift=0cm] 

\coordinate(u) at (4,9);
\coordinate (v) at (16,9);

\node at (3,9) {$u$};
\node at (17,9) {$v$};

\coordinate(v_{11}) at (-3,4);
\coordinate(v_{12}) at (1,4);
\coordinate(v_{13}) at (1,1);
\node at (-4,4.) {$v_{1,1}$};
\node at (1,.5) {$v_{1,2}$};
\node at (2.3,3.5) {$v_{1,3}$};

\draw [thick, blue](v_{11}) -- (u);
\draw [thick, red](v_{11}) -- (v);
\draw [thick, red](v_{12}) -- (v);
\draw [thick, red](v_{13}) -- (v);
\draw (v_{11}) -- (v_{12});
\draw (v_{11}) -- (v_{13});
\draw (v_{12}) -- (v_{13});

\coordinate(v_{21}) at (12,4);
\coordinate(v_{22}) at (12,1);
\coordinate(v_{23}) at (8,4);

\draw [thick, blue](v_{21}) -- (u);
\draw [thick, red](v_{21}) -- (v);
\draw [thick, red](v_{22}) -- (v);
\draw [thick, red](v_{23}) -- (v);

\draw (v_{21}) -- (v_{22});
\draw (v_{21}) -- (v_{23});
\draw (v_{22}) -- (v_{23});

\node at (16,3) {$\cdots$};

\coordinate(v_{k-1 1}) at (19,4);
\coordinate(v_{k-1 2}) at (19,1);
\coordinate(v_{k-1 3}) at (23,4);

\node at (20,4.5) {$v_{k-1, 1}$};
\node at (19.5,.3) {$v_{k-1, 2}$};
\node at (24,4.3) {$v_{k-1, 3}$};

\draw [thick, blue](v_{k-1 1}) -- (u);
\draw [thick, red](v_{k-1 1}) -- (v);
\draw [thick, red](v_{k-1 2}) -- (v);
\draw [thick, red](v_{k-1 3}) -- (v);
\draw (v_{k-1 1}) -- (v_{k-1 2});
\draw (v_{k-1 1}) -- (v_{k-1 3});
\draw (v_{k-1 2}) -- (v_{k-1 3});

\coordinate(v_{k 1}) at (-3,14);
\coordinate(v_{k 2}) at (1,14);
\coordinate(v_{k 3}) at (1,17);
\node at (-4,14) {$v_{k, 1}$};
\node at (0,14.4) {$v_{k, 2}$};
\node at (1,17.9) {$v_{k, 3}$};

\draw [thick, blue](v_{k 1}) -- (u);
\draw [thick, blue](v_{k 3}) -- (u);
\draw [thick, blue](v_{k 2}) -- (u);
\draw [thick, red](v_{k 1}) -- (v);
\draw [thick, red](v_{k 2}) -- (v);

\draw (v_{k 1}) -- (v_{k 2});
\draw (v_{k 1}) -- (v_{k 3});
\draw (v_{k 2}) -- (v_{k 3});

\coordinate(v_{(k+1) 1}) at (7,14);
\coordinate(v_{(k+1) 2}) at (7,17);
\coordinate(v_{(k+1) 3}) at (11,14);

\draw [thick, blue](v_{(k+1) 1}) -- (u);
\draw [thick, blue](v_{(k+1) 3}) -- (u);
\draw [thick, blue](v_{(k+1) 2}) -- (u);
\draw [thick, red](v_{(k+1) 1}) -- (v);
\draw [thick, red](v_{(k+1) 3}) -- (v);

\draw (v_{(k+1) 1}) -- (v_{(k+1) 2});
\draw (v_{(k+1) 1}) -- (v_{(k+1) 3});
\draw (v_{(k+1) 2}) -- (v_{(k+1) 3});

\node at (13,16) {$\cdots$};

\coordinate(v_{(k+l-2) 1}) at (19,14);
\coordinate(v_{(k+l-2) 2}) at (19,17);
\coordinate(v_{(k+l-2) 3}) at (23,14);

\node at (17.5,14.5) {$v_{k+\ell-2, 1}$};
\node at (25,14) {$v_{k+\ell-2, 2}$};
\node at (19,17.7) {$v_{k+\ell-2, 3}$};

\draw [thick, blue](v_{(k+l-2) 1}) -- (u);
\draw [thick, blue](v_{(k+l-2) 3}) -- (u);
\draw [thick, blue](v_{(k+l-2) 2}) -- (u);
\draw [thick, red](v_{(k+l-2) 1}) -- (v);
\draw [thick, red](v_{(k+l-2) 3}) -- (v);
\draw (v_{(k+l-2) 1}) -- (v_{(k+l-2) 2});
\draw (v_{(k+l-2) 1}) -- (v_{(k+l-2) 3});
\draw (v_{(k+l-2) 2}) -- (v_{(k+l-2) 3});
\draw(u)[fill=white] circle(\vr);
\draw(v)[fill=white] circle(\vr);

\draw(v_{11})[fill=white] circle(\vr);
\draw(v_{12})[fill=white] circle(\vr);
\draw(v_{13})[fill=white] circle(\vr);

\draw(v_{21})[fill=white] circle(\vr);
\draw(v_{22})[fill=white] circle(\vr);
\draw(v_{23})[fill=white] circle(\vr);

\draw(v_{k-1 1})[fill=white] circle(\vr);
\draw(v_{k-1 2})[fill=white] circle(\vr);
\draw(v_{k-1 3})[fill=white] circle(\vr);

\draw(v_{k 1})[fill=white] circle(\vr);
\draw(v_{k 2})[fill=white] circle(\vr);
\draw(v_{k 3})[fill=white] circle(\vr);

\draw(v_{(k+1) 1})[fill=white] circle(\vr);
\draw(v_{(k+1) 2})[fill=white] circle(\vr);
\draw(v_{(k+1) 3})[fill=white] circle(\vr);

\draw(v_{(k+l-2) 1})[fill=white] circle(\vr);
\draw(v_{(k+l-2) 2})[fill=white] circle(\vr);
\draw(v_{(k+l-2) 3})[fill=white] circle(\vr);
\end{scope}

\end{tikzpicture}
     \caption{Graph $G_{k,\ell}$ with $\gmb(G_{k,\ell})=k$ and $\gmbt(G_{k,\ell})=\ell$}
    \label{fig:real graph 2}
    \end{center}
\end{figure}

Consider an MBD-game on $G_{k,\ell}$. First, we prove that Dominator has a strategy to finish the D-game on $G_{k,\ell}$ by at most $k$ moves. Dominator selects the vertex $u$ as his first move.  If Staller selects a vertex other than $v$, then Dominator can finish the game by selecting $v$. Thus, Staller selects $v$ as her optimal response. Then the only undominated vertices are in the first $k-1$  triangles. So Dominator can select one vertex from each of the $k-1$ triangles. Hence $\gmb(G_{k,\ell})\leq 1+k-1=k$. Second, we prove that there exists a strategy for Staller which forces Dominator to play at least $k$  moves. If Dominator selects a vertex other than $u$ as his first move, then Staller selects $u$ as her first response. In this case, Dominator needs at least $\ell-1$ more moves to dominate all the vertices of $G_{k,\ell}$, and we have $k\leq \ell$. And if Dominator selects $u$ as his first move,  then Staller selects $v$ as her response. So, Dominator must select at least one vertex from each of the first $k-1$ triangles. Hence, $\gmb(G_{k,\ell})\geq 1+k-1=k$. Therefore, $\gmb(G_{k,\ell})=k$.

We next verify that $\gmbt(G_{k,\ell})=\ell$. Consider the following strategy of Dominator. He selects $v$ as his first move, and then the only undominated vertices are $v_{k,3}, v_{k+1,3},\ldots, v_{k+\ell-2,3}$. Therefore, Dominator can finish the game by selecting either $v_{i,1}$ or $v_{i,2}$ for $i\in\{k,k+1, k+\ell-2\}$. Clearly, the above set of vertices selected by Dominator becomes a total dominating set. Thus, $\gmbt(G_{k,\ell})\leq \ell$. On the other hand, if Dominator first selects a vertex other than $u$, then Staller selects $u$ as her reply. Since any total dominating set without $u$ contains at least $\ell$ vertices, we have $\gmbt(G_{k,\ell}) \ge \ell$. And if Dominator selects $u$ as his first move, then Staller replies by the move $v$. If Dominator selects a vertex in the $j^{\rm th}$ triangle for some $j\in[k+\ell-2]$ as his second move, then Staller selects an unplayed vertex $v_{i,1}$ for $i\neq j$ and $i\in [k-1]$. Note that this is possible because $k\geq 3$. But then Staller wins the game by selecting either $v_{i,2}$ or $v_{i,3}$ after the next move of Dominator. We can conclude that $\gmbt(G_{k,\ell})\geq \ell$ and therefore $\gmbt(G_{k,\ell})=\ell$.
\qed

The next result is the analogous result to Theorem~\ref{thm:MB-MBT} for the S-game. 

\begin{theorem}
\label{thm:MB'-MBT'}
For any two integers $k$ and $\ell$ with $2\leq k\leq  \ell$, there exists a connected graph $G$ with $\gmb'(G)=k$ and $\gmbt'(G)=\ell$.
\end{theorem}

\proof
Clearly, $\gmb'(C_4)=2$  and $\gmbt'(C_4)=2$ and this proves the result  for $\ell=2$. It remains to consider case when $\ell\geq 3$. 

\begin{figure}[ht!]
\begin{center}
\begin{tikzpicture}[scale=0.45,style=thick,x=1cm,y=1cm]
\def\vr{6pt}
\begin{scope}[xshift=0cm, yshift=0cm] 

\coordinate(u) at (3,9);
\coordinate(v) at (9,9);
\coordinate (w) at (16,9);

\node at (2,9) {$u$};
\node at (9.5,8.5) {$v$};
\node at (17,9) {$w$};

\draw [thick, red](u) -- (v);
\draw [thick, green] (u) to[out=-20,in=-70] (w);

\coordinate(v_{12}) at (11,1);
\coordinate(v_{11}) at (8,4);
\coordinate(v_{13}) at (5.5,1);

\node at (4,.9) {$v_{1, 3}$};
\node at (12,1) {$v_{1, 2}$};
\node at (9.3,3.7) {$v_{1, 1}$};

\draw[thick, blue](v_{11}) -- (u);
\draw [thick, red](v_{11}) -- (v);
\draw [thick, red](v_{12}) -- (v);
\draw [thick, red](v_{13}) -- (v);
\draw [thick, green](v_{11}) -- (w);

\draw (v_{11}) -- (v_{12});
\draw (v_{11}) -- (v_{13});
\draw (v_{12}) -- (v_{13});

\coordinate(v_{2 1}) at (1,17);
\coordinate(v_{2 2}) at (1,14);
\coordinate(v_{2 3}) at (-4,14);

\node at (.3,14.5) {$v_{2, 1}$};
\node at (1,17.5) {$v_{2, 2}$};
\node at (-5,14) {$v_{2, 3}$};

\draw [thick, blue](v_{2 1}) -- (u);
\draw [thick, blue](v_{2 3}) -- (u);
\draw [thick, blue](v_{2 2}) -- (u);
\draw [thick, red](v_{2 1}) -- (v);
\draw [thick, red](v_{2 3}) -- (v);
\draw [thick, red](v_{2 2}) -- (v);
\draw [thick, green](v_{2 1}) -- (w);
\draw [thick, green](v_{2 2}) -- (w);

\draw (v_{2 1}) -- (v_{2 2});
\draw (v_{2 1}) -- (v_{2 3});
\draw (v_{2 2}) -- (v_{2 3});

\coordinate(v_{3 1}) at (9.5,14);
\coordinate(v_{3 2}) at (12,17);
\coordinate(v_{3 3}) at (7,17);

\draw [thick, blue](v_{3 1}) -- (u);
\draw [thick, blue](v_{3 3}) -- (u);
\draw [thick, blue](v_{3 2}) -- (u);
\draw [thick, red](v_{3 1}) -- (v);
\draw [thick, red](v_{3 3}) -- (v);
\draw [thick, red](v_{3 2}) -- (v);
\draw [thick, green](v_{3 1}) -- (w);
\draw [thick, green](v_{3 2}) -- (w);

\draw (v_{3 1}) -- (v_{3 2});
\draw (v_{3 1}) -- (v_{3 3});
\draw (v_{3 2}) -- (v_{3 3});

\node at (14,15) {$\cdots$};
\node at (14,17) {$\cdots$};

\coordinate(v_{l-1 1}) at (19,14);
\coordinate(v_{l-1 2}) at (23,14);
\coordinate(v_{l-1 3}) at (19,17);
\node at (17.5,14.4) {$v_{\ell-1, 1}$};
\node at (23.9,14.7) {$v_{\ell-1, 2}$};
\node at (19,17.8) {$v_{\ell-1, 3}$};

\draw [thick, blue](v_{l-1 1}) -- (u);
\draw [thick, blue](v_{l-1 3}) -- (u);
\draw [thick, blue](v_{l-1 2}) -- (u);
\draw [thick, red](v_{l-1 1}) -- (v);
\draw [thick, red](v_{l-1 3}) -- (v);
\draw [thick, red](v_{l-1 2}) -- (v);
\draw [thick, green](v_{l-1 1}) -- (w);
\draw [thick, green](v_{l-1 2}) -- (w);

\draw (v_{l-1 1}) -- (v_{l-1 2});
\draw (v_{l-1 1}) -- (v_{l-1 3});
\draw (v_{l-1 2}) -- (v_{l-1 3});
\draw(u)[fill=white] circle(\vr);
\draw(v)[fill=white] circle(\vr);
\draw(w)[fill=white] circle(\vr);

\draw(v_{11})[fill=white] circle(\vr);
\draw(v_{12})[fill=white] circle(\vr);
\draw(v_{13})[fill=white] circle(\vr);

\draw(v_{2 1})[fill=white] circle(\vr);
\draw(v_{2 2})[fill=white] circle(\vr);
\draw(v_{2 3})[fill=white] circle(\vr);

\draw(v_{3 1})[fill=white] circle(\vr);
\draw(v_{3 2})[fill=white] circle(\vr);
\draw(v_{3 3})[fill=white] circle(\vr);

\draw(v_{l-1 1})[fill=white] circle(\vr);
\draw(v_{l-1 2})[fill=white] circle(\vr);
\draw(v_{l-1 3})[fill=white] circle(\vr);

\end{scope}

\end{tikzpicture}
     \caption{Graph $H_{2, \ell}$ with $\gmb'(H_{2, \ell})=2$ and $\gmbt'(H_{2, \ell})=\ell$}
    \label{fig:real graph 3}
    \end{center}
\end{figure}
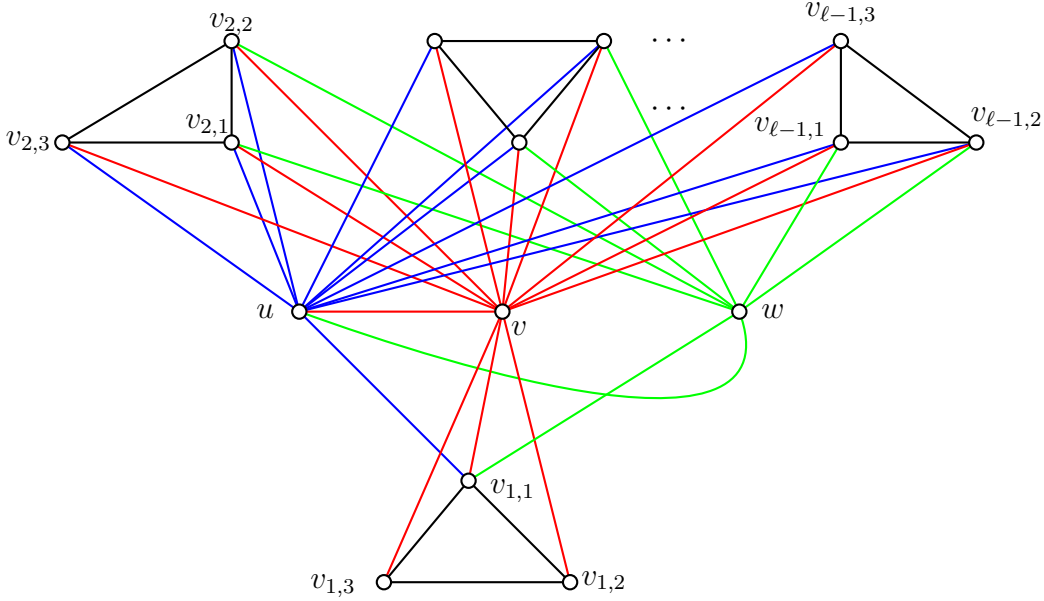

\medskip\noindent
{\bf Case 1}: $k = 2$. \\
Let $H_{2,\ell}$, $\ell\ge 3$, be the graph constructed in the following way. First, take $\ell-1$  disjoint triangles with the vertices $v_{i,j}$, $j\in [3]$, in the $i^{\rm th}$ triangle, $i\in[\ell-1]$. Next, take three additional vertices $u$, $v$, and $w$. The vertex $v$ is adjacent to all the other vertices but $w$. The vertex $u$ is adjacent to all the other vertices but $v_{1,2}$, and $v_{1,3}$. Finally, $w$ is adjacent to the vertices $v_{i,j}$, $i\in \{2, 3, \ldots, \ell-1\}$, $j\in [2]$, and to the vertices $u$ and $v_{1,1}$. See Fig.~\ref{fig:real graph 3}.

We first claim that $\gamma'_{MB}(H_{2, \ell})=2$. Indeed, Dominator can select one of $u$ and $v$ as his first move. If he selects $u$, then he finishes the game by selecting one vertex from the set $\{ v_{1,1}, v_{1,2}, v_{1,3}\}$ as his second move. If he first selects $v$, then he wins in his second move by selecting either $w$ or a vertex $v_{i,1}$ for some $i\in[\ell-1]$. Therefore, $\gmb'(H_{2, \ell})=2$.

Consider now an MBTD game played on $H_{2, \ell}$. First,  we show that $\gmbt'(H_{2, \ell})\leq \ell$. If Staller selects a vertex other than $v$, Then Dominator can finish the game by selecting the vertex $v$ and any unplayed vertex adjacent to $w$. Therefore, Staller prefers to select $v$ as her first move. This forces Dominator to select $v_{1,1}$ as his response, for otherwise Staller can win on the first triangle. Next Staller selects $u$ as her second, for otherwise Dominator selects $u$ as his second move and then finishes the game.  Then Dominator selects $w$, and later on he wins by selecting one of $v_{i,1}$ and $v_{i,2}$ for each $i\in\{2, 3, \ldots, \ell-1\}$. Therefore, he has played at most $2+\ell-2=\ell$ moves, that is, $\gmb'(H_{2, \ell})\leq \ell$.

Consider next the Staller's view point. Let she starts the game by selecting $v$. If Dominator selects any vertex other than $v_{1,1}$ as his reply, then Staller plays on $v_{1,1}$, and she wins the game by selecting an unplayed vertex from the set $\{v_{1,2}, v_{1,3}\}$ in one more move. So, Dominator is forced to select $v_{1,1}$ as his first move. Then Staller selects $u$ as her second move. If Dominator selects any vertex other than $w$, then Staller replies by selecting the vertex  $w$. Note that every total dominating set without $u, v$,  and $w$ contains at least $2(\ell-1)$ vertices. Since $\ell\geq 3$, we have $2(\ell-1)\geq \ell$. Therefore, Dominator selects the vertex $w$ as his next move. It is now evident that Dominator requires to select either $v_{i,1}$ or $v_{i,2}$ for $i \in\{ 2,3,  \ldots, \ell-1\}$ to win the game. Therefore, $\gmbt'(H_{2, \ell})\geq 2+\ell-2 =\ell$. Thus, $\gmbt'(H_{2, \ell}) = \ell$.

\medskip\noindent
{\bf Case 2}: $k\ge 3$. \\
Let $H_{k,\ell}$, $\ell\ge k\ge 3$, be the graph constructed in the following way. First take $k+\ell-2$ disjoint triangles with vertices $v_{i,j}$, $j\in [3]$, in the $i^{\rm th}$ triangle, $i\in [k+\ell - 2]$. Next, take three more vertices $u$, $v$, and $w$. The vertex $u$ is adjacent to $v$, to the vertices $v_{i,1}$, $i\in [k+\ell-2]$, and to the vertices $v_{i,j}$, $i\in \{k, k+1, \ldots, k+\ell-2\}$, $j\in \{2, 3\}$. The vertex $v$ is adjacent to $u$, to the vertices $v_{i,j}$, $i\in [k+\ell-2]$, $j\in [2]$, and to the vertices $v_{i,3}$, $i\in \{k, k+1, \ldots, k+\ell-2\}$. Finally, the vertex $w$ is adjacent to the vertices $v_{i,j}$, $i\in [k+\ell-2]$, $j\in [2]$, and to the vertices $v_{i,3}$, $i\in [k-1]$. See Fig.~\ref{fig:real graph 4}. 

\begin{figure}[ht!]
\begin{center}
\begin{tikzpicture}[scale=0.45,style=thick,x=1cm,y=1cm]
\def\vr{6pt}
\begin{scope}[xshift=0cm, yshift=0cm] 

\coordinate(u) at (3,9);
\coordinate(v) at (8,9);
\coordinate (w) at (16,9);

\node at (9.5,9) {$u$};
\node at (2,9) {$v$};
\node at (17,9) {$w$};
\draw [thick, blue](u) -- (v);

\coordinate(v_{11}) at (1,4);
\coordinate(v_{12}) at (-3,4);
\coordinate(v_{13}) at (1,1);
\node at (2,3.5) {$v_{1,1}$};
\node at (-4,4.5) {$v_{1,2}$};
\node at (1,0) {$v_{1,3}$};

\draw [thick, blue](v_{11}) -- (u);
\draw [thick, blue](v_{12}) -- (u);
\draw [thick, red](v_{11}) -- (v);
\draw [thick, green](v_{11}) -- (w);
\draw [thick, green](v_{12}) -- (w);
\draw [thick, green](v_{13}) -- (w);

\draw (v_{11}) -- (v_{12});
\draw (v_{11}) -- (v_{13});
\draw (v_{12}) -- (v_{13});

\coordinate(v_{21}) at (12,4);
\coordinate(v_{22}) at (8,4);
\coordinate(v_{23}) at (12,1);

\draw [thick, blue](v_{21}) -- (u);
\draw [thick, blue](v_{22}) -- (u);
\draw [thick, red](v_{21}) -- (v);
\draw [thick, green](v_{21}) -- (w);
\draw [thick, green](v_{22}) -- (w);
\draw [thick, green](v_{23}) -- (w);

\draw (v_{21}) -- (v_{22});
\draw (v_{21}) -- (v_{23});
\draw (v_{22}) -- (v_{23});

\node at (15.5,2.5) {$\cdots$};

\coordinate(v_{k-1 1}) at (19,4);
\coordinate(v_{k-1 2 }) at (19,1);
\coordinate(v_{k-1 3}) at (23,4);
\node at (17.5,3.5) {$v_{k-1, 1}$};
\node at (24.5,4) {$v_{k-1, 2}$};
\node at (20.5,.5) {$v_{k-1, 3}$};

\draw [thick, blue](v_{k-1 1}) -- (u);
\draw [thick, blue](v_{k-1 3}) -- (u);
\draw [thick, red](v_{k-1 1}) -- (v);
\draw [thick, green](v_{k-1 1}) -- (w);
\draw [thick, green](v_{k-1 2}) -- (w);
\draw [thick, green](v_{k-1 3}) -- (w);

\draw (v_{k-1 1}) -- (v_{k-1 2});
\draw (v_{k-1 1}) -- (v_{k-1 3});
\draw (v_{k-1 2}) -- (v_{k-1 3});

\coordinate(v_{k 1}) at (1,17);
\coordinate(v_{k 2}) at (1,14);
\coordinate(v_{k 3}) at (-3,14);

\node at (0,14.5) {$v_{k, 1}$};
\node at (1,18) {$v_{k, 2}$};
\node at (-4,14.5) {$v_{k, 3}$};

\draw [thick, blue](v_{k 1}) -- (u);
\draw [thick, blue](v_{k 3}) -- (u);
\draw [thick, blue](v_{k 2}) -- (u);
\draw [thick, red](v_{k 1}) -- (v);
\draw [thick, red](v_{k 3}) -- (v);
\draw [thick, red](v_{k 2}) -- (v);
\draw [thick, green](v_{k 1}) -- (w);
\draw [thick, green](v_{k 2}) -- (w);

\draw (v_{k 1}) -- (v_{k 2});
\draw (v_{k 1}) -- (v_{k 3});
\draw (v_{k 2}) -- (v_{k 3});

\coordinate(v_{k+1 1}) at (9.5,14);
\coordinate(v_{k+1 2}) at (12,17);
\coordinate(v_{k+1 3}) at (7,17);

\draw [thick, blue](v_{k+1 1}) -- (u);
\draw [thick, blue](v_{k+1 3}) -- (u);
\draw [thick, blue](v_{k+1 2}) -- (u);
\draw [thick, red](v_{k+1 1}) -- (v);
\draw [thick, red](v_{k+1 3}) -- (v);
\draw [thick, red](v_{k+1 2}) -- (v);
\draw [thick, green](v_{k+1 1}) -- (w);
\draw [thick, green](v_{k+1 2}) -- (w);

\draw (v_{k+1 1}) -- (v_{k+1 2});
\draw (v_{k+1 1}) -- (v_{k+1 3});
\draw (v_{k+1 2}) -- (v_{k+1 3});

\node at (15,16.5) {$\cdots$};

\coordinate(v_{k+l-2 1}) at (19,14);
\coordinate(v_{k+l-2 2}) at (23,14);
\coordinate(v_{k+l-2 3}) at (19,17);
\node at (17.3,14.4) {$v_{k+\ell-2, 1}$};
\node at (24.5,14.4) {$v_{k+\ell-2, 2}$};
\node at (19,18) {$v_{k+\ell-2, 3}$};

\draw [thick, blue](v_{k+l-2 1}) -- (u);
\draw [thick, blue](v_{k+l-2 3}) -- (u);
\draw [thick, blue](v_{k+l-2 2}) -- (u);
\draw [thick, red](v_{k+l-2 1}) -- (v);
\draw [thick, red](v_{k+l-2 3}) -- (v);
\draw [thick, red](v_{k+l-2 2}) -- (v);
\draw [thick, green](v_{k+l-2 1}) -- (w);
\draw [thick, green](v_{k+l-2 2}) -- (w);

\draw (v_{k+l-2 1}) -- (v_{k+l-2 2});
\draw (v_{k+l-2 1}) -- (v_{k+l-2 3});
\draw (v_{k+l-2 2}) -- (v_{k+l-2 3});
\draw(u)[fill=white] circle(\vr);
\draw(v)[fill=white] circle(\vr);
\draw(w)[fill=white] circle(\vr);

\draw(v_{11})[fill=white] circle(\vr);
\draw(v_{12})[fill=white] circle(\vr);
\draw(v_{13})[fill=white] circle(\vr);

\draw(v_{k-1 1})[fill=white] circle(\vr);
\draw(v_{k-1 2})[fill=white] circle(\vr);
\draw(v_{k-1 3})[fill=white] circle(\vr);

\draw(v_{21})[fill=white] circle(\vr);
\draw(v_{22})[fill=white] circle(\vr);
\draw(v_{23})[fill=white] circle(\vr);

\draw(v_{k 1})[fill=white] circle(\vr);
\draw(v_{k 2})[fill=white] circle(\vr);
\draw(v_{k 3})[fill=white] circle(\vr);

\draw(v_{k+1 1})[fill=white] circle(\vr);
\draw(v_{k+1 2})[fill=white] circle(\vr);
\draw(v_{k+1 3})[fill=white] circle(\vr);

\draw(v_{k+l-2 1})[fill=white] circle(\vr);
\draw(v_{k+l-2 2})[fill=white] circle(\vr);
\draw(v_{k+l-2 3})[fill=white] circle(\vr);

\end{scope}

\end{tikzpicture}
     \caption{Graph $H_{k,\ell}$ with $\gmb'(H_{k,\ell})=k$ and $\gmbt'(H_{k,\ell})=\ell$}
    \label{fig:real graph 4}
    \end{center}
\end{figure}

Dominator can select one of $u$ and $v$ as his first move, thus dominating all the vertices in the $i^{\rm th}$ triangle for $i\in\{k, k+1, \ldots, k+\ell-2\}$ is dominated. So Dominator can finish the game by selecting one of $v_{i,1}$ and $v_{i,2}$ for $i \in [k-1]$. Hence $\gmb'(H_{k,\ell})\leq 1+k-1=k$. To show that $\gmb'(H_{k,\ell})\geq k$, consider the following strategy of Staller. She selects $v$ as her first move. Since every dominating set of $H_{k,\ell}$ without $u$, $v$, and $w$ contains at least $k+\ell-2$ vertices, Dominator replies by selecting $u$ or $w$. If he selects  $u$, then Staller's second move is $w$, while if Dominator's first move is $w$, Staller's second move is $u$. In the first case Dominator needs at least $k-1$ more moves to finish the game, and in the second case at least $\ell-1\geq k-1$ more moves. Therefore $\gmb'(H_{k,\ell})\geq k$ and hence $\gmb'(H_{k,\ell})= k$.

It remains to prove that $\gmbt'(H_{k,\ell}) = \ell$. First, consider the following strategy of  Dominator. If Staller start the game by selecting $w$, Dominator replies by $v$. Then Dominator can win the game by selecting one of $v_{i,1}$ and $v_{i,2}$  for $i\in [k-1]$, so he needs at most $1+k-1=k\leq \ell$ moves. If Staller selects the vertex $u$ or $v$ as her first move, then Dominator selects $w$ as his response. Afterwards,  Dominator can win the game by selecting one of $v_{i,1}$ and $v_{i,2}$  for every $i\in\{k, k+1, \ldots, k+\ell-2\}$. Therefore, he can finish the game in at most $1+\ell-1 = \ell$ moves. The same argument applies when Staller selects a vertex $v_{i,j}$,  where $i\in [k-1]$ and $j\in [3]$, as her first move. Finally assume that Staller first selects $v_{i,j}$, where $k\leq i\leq k+\ell-2$ and $1\leq j\leq 3$. In this case  Dominator selects $v$ as his reply. After that, he can win by selecting one of $v_{i,1}$ and $v_{i,2}$  for every $i\in[k-1]$, that is, after at most $1+k-1 = k\leq \ell$ moves. We can conclude that $\gmbt'(H_{k,\ell})\leq \ell$.

To show that $\gmbt'(H_{k,\ell})\geq \ell$, assume that Staller first selects $v$. If Dominator replies by $u$, then Staller's second move is $w$. Now Dominator plays a vertex $v_{i,1}$, where $i\in [k-1]$, and Staller selects a vertex $v_{j,1}$, where $j\in [k-1]$ and $j\neq i$. Then Staller can win the game by selecting one more vertex from the $j^{\rm th}$ triangle. Note that a total dominating set which does not contain $u$, $v$, and $w$, contains at least $2(k+\ell-2)$ vertices. Therefore, Dominator selects $w$ as his first move. In that case, Staller selects $u$ as her second move. But then Dominator needs to select one of $v_{i,1}$ and $v_{i,2}$ for $i\in\{k, k+1, \ldots, k+\ell-2\}$. We can conclude that $\gmbt'(H_{k,\ell})\geq \ell$ and thus $\gmbt'(H_{k,\ell}) = \ell$.
\qed 

Our final realization result reads as follows. 

\begin{theorem}\label{thm:MBT-MBT'}
For any two integers $k$ and $\ell$ with $2\leq k\leq \ell$, there exists a connected graph $G$ with $\gmbt(G)=k$ and $\gmbt'(G)=\ell$. 
\end{theorem}

\begin{figure}[ht!]
\begin{center}
\begin{tikzpicture}[scale=0.45,style=thick,x=1cm,y=1cm]
\def\vr{6pt}
\begin{scope}[xshift=0cm, yshift=0cm] 
\coordinate(u) at (4,9);
\coordinate (v) at (16,9);
\node at (3,9) {$u$};
\node at (17,9) {$v$};

\coordinate(v_{11}) at (-4,4);
\coordinate(v_{12}) at (1,4);
\coordinate(v_{13}) at (1,0);
\node at (2, 3.5) {$v_{1,1}$};
\node at (-4.5,4.5) {$v_{1,2}$};
\node at (.5,-.5) {$v_{1,3}$};

\draw [thick, blue](v_{11}) -- (u);
\draw [thick, blue](v_{12}) -- (u);
\draw [thick, red](v_{11}) -- (v);
\draw [thick, red](v_{12}) -- (v);
\draw [thick, red](v_{13}) -- (v);
\draw (v_{11}) -- (v_{12});
\draw (v_{11}) -- (v_{13});
\draw (v_{12}) -- (v_{13});

\coordinate(v_{21}) at (7,0);
\coordinate(v_{22}) at (12,0);
\coordinate(v_{23}) at (12,4);
\node at (13, 4) {$v_{2,1}$};
\node at (13,0) {$v_{2,3}$};
\node at (6,-.5) {$v_{2,2}$};

\draw [thick, blue](v_{23}) -- (u);
\draw [thick, blue](v_{21}) -- (u);
\draw [thick, red](v_{21}) -- (v);
\draw [thick, red](v_{22}) -- (v);
\draw [thick, red](v_{23}) -- (v);

\draw (v_{21}) -- (v_{22});
\draw (v_{21}) -- (v_{23});
\draw (v_{22}) -- (v_{23});

\node at (15.5,2) {$\cdots$};

\coordinate(v_{k-1 1}) at (19,4);
\coordinate(v_{k-1 2}) at (24,4);
\coordinate(v_{k-1 3}) at (19,0);
\node at (20.5,3.5) {$v_{k-1, 1}$};
\node at (25.5,4) {$v_{k-1, 2}$};
\node at (19,-.5) {$v_{k-1, 3}$};

\draw [thick, blue](v_{k-1 2}) -- (u);
\draw [thick, blue](v_{k-1 1}) -- (u);
\draw [thick, red](v_{k-1 1}) -- (v);
\draw [thick, red](v_{k-1 2}) -- (v);
\draw [thick, red](v_{k-1 3}) -- (v);

\draw (v_{k-1 1}) -- (v_{k-1 2});
\draw (v_{k-1 1}) -- (v_{k-1 3});
\draw (v_{k-1 2}) -- (v_{k-1 3});

\coordinate(v_{k 1}) at (1,14);
\coordinate(v_{k 2}) at (1,18);
\coordinate(v_{k 3}) at (-4,14);
\node at (.2,14.5) {$v_{k, 1}$};
\node at (1,19) {$v_{k, 2}$};
\node at (-5.5,14.5) {$v_{k, 3}$};

\draw [thick, blue](v_{k 1}) -- (u);
\draw [thick, blue](v_{k 2}) -- (u);
\draw [thick, blue](v_{k 3}) -- (u);
\draw [thick, red](v_{k 1}) -- (v);
\draw [thick, red](v_{k 2}) -- (v);

\draw (v_{k 1}) -- (v_{k 2});
\draw (v_{k 1}) -- (v_{k 3});
\draw (v_{k 2}) -- (v_{k 3});

\coordinate(u_{(k+1)1}) at (7.5,15);
\coordinate(u_{(k+1)2}) at (7.5,19);
\coordinate(u_{(k+1)3}) at (12.5,15);

\draw [thick, blue](u_{(k+1)1}) -- (u);
\draw [thick, blue](u_{(k+1)3}) -- (u);
\draw [thick, blue](u_{(k+1)2}) -- (u);
\draw [thick, red](u_{(k+1)1}) -- (v);
\draw [thick, red](u_{(k+1)3}) -- (v);

\draw (u_{(k+1)1}) -- (u_{(k+1)2});
\draw (u_{(k+1)1}) -- (u_{(k+1)3});
\draw (u_{(k+1)2}) -- (u_{(k+1)3});

\node at (14,18) {$\cdots$};

\coordinate(u_{k+l-2 1}) at (19,15);
\coordinate(u_{k+l-2 2}) at (19,19);
\coordinate(u_{k+l-2 3}) at (24,15);
\node at (21,15.5) {$v_{k+\ell-2, 1}$};
\node at (19.5,19.5) {$v_{k+\ell-2, 2}$};
\node at (25.5,15.5) {$v_{k+\ell-2, 3}$};

\draw [thick, blue](u_{k+l-2 1}) -- (u);
\draw [thick, blue](u_{k+l-2 3}) -- (u);
\draw [thick, blue](u_{k+l-2 2}) -- (u);
\draw [thick, red](u_{k+l-2 1}) -- (v);
\draw [thick, red](u_{k+l-2 2}) -- (v);
\draw (u_{k+l-2 1}) -- (u_{k+l-2 2});
\draw (u_{k+l-2 1}) -- (u_{k+l-2 3});
\draw (u_{k+l-2 2}) -- (u_{k+l-2 3});
\draw(u)[fill=white] circle(\vr);
\draw(v)[fill=white] circle(\vr);

\draw(v_{11})[fill=white] circle(\vr);
\draw(v_{12})[fill=white] circle(\vr);
\draw(v_{13})[fill=white] circle(\vr);

\draw(v_{21})[fill=white] circle(\vr);
\draw(v_{22})[fill=white] circle(\vr);
\draw(v_{23})[fill=white] circle(\vr);

\draw(v_{k 1})[fill=white] circle(\vr);
\draw(v_{k 2})[fill=white] circle(\vr);
\draw(v_{k 3})[fill=white] circle(\vr);

\draw(v_{k-1 1})[fill=white] circle(\vr);
\draw(v_{k-1 2})[fill=white] circle(\vr);
\draw(v_{k-1 3})[fill=white] circle(\vr);

\draw(u_{(k+1)1})[fill=white] circle(\vr);
\draw(u_{(k+1)2})[fill=white] circle(\vr);
\draw(u_{(k+1)3})[fill=white] circle(\vr);

\draw(u_{k+l-2 1})[fill=white] circle(\vr);
\draw(u_{k+l-2 2})[fill=white] circle(\vr);
\draw(u_{k+l-2 3})[fill=white] circle(\vr);
\end{scope}

\end{tikzpicture}
     \caption{Graph $F_{k,\ell}$ with $\gmbt(F_{k,\ell})=k$ and $\gmbt'(F_{k,\ell})=\ell$}
    \label{fig: real graph 5}
    \end{center}
\end{figure}
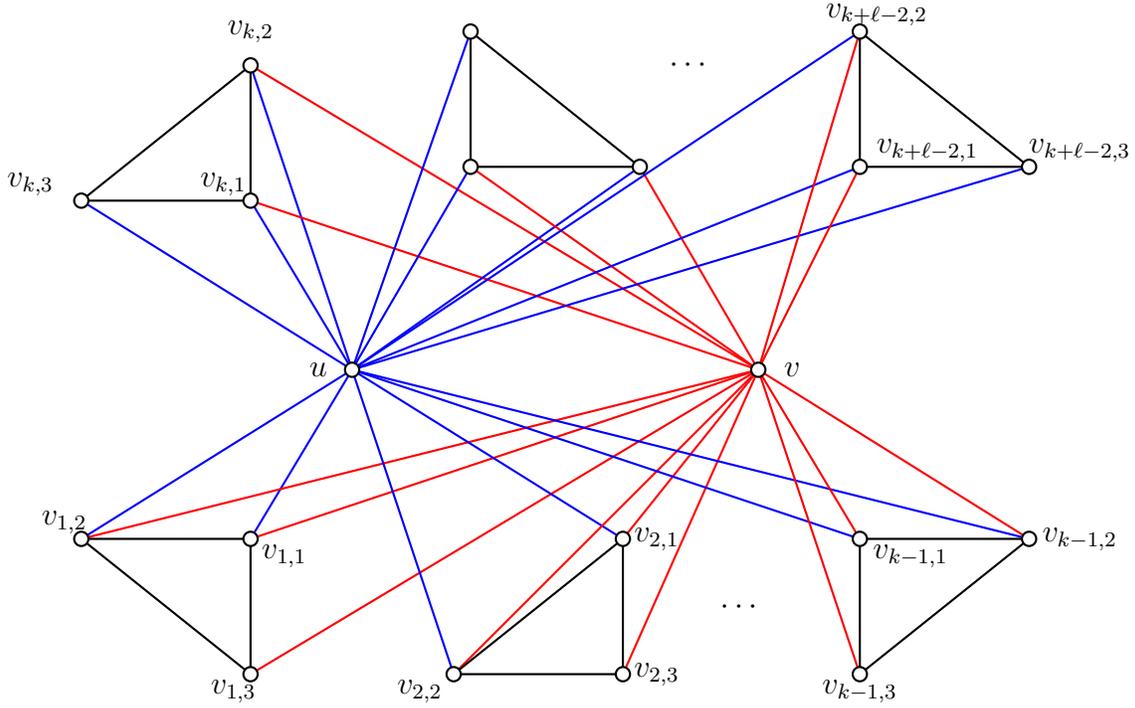

\proof
Let $F_{k,\ell}$, $\ell\ge k\ge 2$, be the graph constructed in the following way. First, take $k+\ell-2$ disjoint triangles, where $v_{i,j}$, $j\in [3]$, are the vertices of the $i^{\rm th}$ triangle, $i\in[k+\ell-2]$. 
Next, take two more vertices $u$ and $v$. The vertex $u$ is adjacent to $v_{i,j}$, $i\in [k+\ell-2]$, $j\in [2]$, and to the vertices $v_{i,3}$, $i\in \{k, k+1,\ldots, k+\ell-2\}$. The vertex $v$ is adjacent to $v_{i,j}$,  $i\in [k+\ell-2]$, $j\in [2]$, and to $v_{i,3}$, $i\in [k-1]$. See Fig.~\ref{fig: real graph 5}. 

To show that $\gmbt(F_{k,\ell})\leq k$, assume that Dominator selects $u$ as his first  move. If Staller selects a vertex other than $v$, then Dominator can finish the game by selecting $v$, and subsequently by playing an unplayed vertex $v_{i,1}$ for some $i$. Therefore, Staller is forced to select $v$ first. Then Dominator can finish the game by selecting one of the vertices from $\{v_{i,1}, v_{i,2}\}$ for $i\in[k-1]$. Hence, Dominator can finish the game in at most $1+k-1=k$ moves, hence $\gmbt(F_{k,\ell})\leq k$.
 
To prove that $\gmbt(F_{k,\ell})\geq k$, observe that any total dominating set of $F_{k,\ell}$ without the vertices $u$ and $v$ contains at least $2(k+\ell-2)$ vertices. Since $k\leq \ell$ and $k\geq 2$, we have $2(k+\ell-2) > k$. Therefore, Dominator starts the game by playing either $v$ or $u$. In the first case Staller's reply is $u$. Then the vertices $v_{i,3}$, $i\in\{k, k+1, \ldots, k+\ell-2\}$, are undominated. Therefore, Dominator needs at least $\ell-1$ more moves to finish the game. And if Dominator selects $u$ as his first move, then Staller selects $v$ as her response. In this case, the vertices $v_{i,3}$, $i\in[k-1]$, remain undominated, and Dominator needs at least $k-1$ more moves to finish the game. Hence $\gmbt(F_{k,\ell})\geq 1 + k-1=k$ which implies that $\gmbt(F_{k,\ell})=k$.

It remains to prove that $\gmbt'(F_{k,\ell})=\ell$. Assume that Staller starts the game by  selecting $u$. Note that any total dominating set without $u$ and $v$ contains at least $2(k+\ell-2)$ vertices. Since $k\geq 2$ and $k\leq \ell$, we have $2(k+\ell-2)>\ell$. Thus, Dominator must select $v$ as his first move, for otherwise he needs at least $2(k+\ell-2)\geq \ell$  moves to obtain a total dominating set. Now, Dominator can win the game by selecting one vertex from $\{v_{i,1}, v_{i,2}\}$ for $i\in\{k, k+1, \ldots, k+\ell-2\}$. Hence Dominator needs at least $\ell-1$ more moves to get a total dominating set. Thus $\gmbt'(F_{k,\ell})\geq \ell$. Finally, consider the following strategy of Dominator. Assume first that Staller selects a vertex from the $i^{\rm th}$ triangle, where $i\in [k-1]$. Then Dominator replies by the move $v$. This enables him to win the game by selecting one of $v_{i,1}$ and $v_{i,2}$ for $i\in\{k, k+1, \ldots, k+l-2\}$. Assume second that Staller first selects a vertex from the $i^{\rm th}$ triangle, where $k\leq i\leq k+\ell-2$. In this case, Dominator selects $u$ as his first move, and he can win the game by selecting one of $v_{i,1}$ and $v_{i,2}$ for $i\in[k-1]$. Thus, the game has at most $1+k-1=k\leq \ell$ moves. Assume finally that Staller selects $u$ or $v$ as her first move. In this case Dominator selects the other of these two vertices in his first move. In both cases he can finish the game in at most $\ell$ moves. Thus $\gmbt'(F_{k,\ell})\leq \ell$ and we can conclude that $\gmbt'(F_{k,\ell})=\ell$.
\qed



\section*{Acknowledgements}

Sand Klav\v{z}ar was supported by the Slovenian Research and Innovation Agency (ARIS) under the grants P1-0297, N1-0285, and N1-0355. This research was partly carried out during the Workshop on Games on Graphs III, June 2025, Rogla, Slovenia, the authors thank the Institute of Mathematics, Physics and Mechanics, Ljubljana, Slovenia for supporting the workshop. 

 
 



\begin{thebibliography}{99}

\bibitem{bagan-2025}
G.~Bagan, E.~Duch\^{e}ne, V.~Gledel, T.~Lehtil\"a, A.~Parreau, 
Partition strategies for the Maker-Breaker domination game, 
Algorithmica 87 (2025) 191--222.

\bibitem{bresar-2010}
B.~Bre{\v{s}}ar, S.~Klav\v{z}ar, D.F.~Rall,
Domination game and an imagination strategy,
SIAM J.\ Discrete Math.\ 24 (2010) 979--991.
	
\bibitem{book-2021}
B.~Bre\v{s}ar, M.A.~Henning, S.~Klav\v{z}ar, D.F.~Rall, 
Domination Games Played on Graphs, 
SpringerBriefs in Mathematics, 2021.  

\bibitem{bujtas-2024} 
Cs.~Bujt\'as, P.~Dokyeesun, 
Fast winning strategies for Staller in the Maker-Breaker domination game, 
Discrete Appl.\ Math.\ 344 (2024) 10--22.

\bibitem{bujtas-2023} 
Cs.~Bujt\'as, P.~Dokyeesun, S.~Klav\v{z}ar, 
Maker-Breaker domination game on trees when Staller wins,  
Discrete Math.\ Theor.\ Comput.\ Sci.\ 25(2) (2023) Paper 12.

\bibitem{divakaran-2025}
A.~Divakaran, T.~Dravec, T.~James, S.~Klav\v{z}ar, L.S.~Nair, 
Maker-Breaker domination game critical graphs,
Discrete Appl.\ Math.\ 368 (2025) 126--134.

\bibitem{dokyeesun-2024} 
P.~Dokyeesun, 
Contributions to Maker-Breaker Domination Game,
Doctoral Thesis, University of Ljubljana, Ljubljana, 2024. 

\bibitem{dokyeesun-2024+} 
P.~Dokyeesun, 
Maker-Breaker domination game on Cartesian products of graphs, 
Commun.\ Comb.\ Optim.\ (2025) \url{doi.org/10.22049/cco.2025.29866.2198}.

\bibitem{duchene-2020}
E.~Duch\^{e}ne, V.~Gledel, A.~Parreau, G.~Renault,
Maker-Breaker domination game,
Discrete Math.\ 343 (2020) 111955.

\bibitem{erdos-1973} 
P.~Erd\H{o}s, J.L.~Selfridge, 
On a combinatorial game, 
J.\ Combin.\ Theory Ser.\ A 14 (1973) 298--301.

\bibitem{favaron-2015}
O.~Favaron, H.~Karami, S.M.~Sheikholeslami, 
Game domination subdivision number of a graph,
J.\ Comb.\ Optim.\ 30 (2015) 109--119. 

\bibitem{forcan-2022}
J.~Forcan, M.~Mikala\v{c}ki,
Maker-Breaker total domination game on cubic graphs,
Discrete Math.\ Theor.\ Comput.\ Sci.\ 24(1) (2022) Paper 20.

\bibitem{forcan-2023}
J.~Forcan, J.~Qi,
Maker-Breaker domination number for Cartesian products of path graphs $P_2$ and $P_n$, 
Discrete Math.\ Theor.\ Comput.\ Sci.\ 25(2) (2023) Paper 21.

\bibitem{gledel-2020} 
V.~Gledel, M.A.~Henning, V. Ir\v{s}i\v{c}, S.~Klav\v{z}ar, 
Maker-Breaker total domination game, 
Discrete Appl.\ Math.\ 282 (2020) 96--107.

\bibitem{gledel-2019} 
V.~Gledel,  V.~Ir\v{s}i\v{c}, S.~Klav\v{z}ar, 
Maker-Breaker domination number, 
Bull.\ Malays.\ Math.\ Sci.\ Soc.\ 42 (2019) 1773--1789.

\bibitem{HJ-1963} 
A.W.\ Hales, R.I.\ Jewett, 
Regularity and positional games, 
Trans.\ Amer.\ Math.\ Soc.\ 106 (1963) 222--229.


\bibitem{hefetz-2014}
D.~Hefetz, M.~Krivelevich, M.~Stojakovi\'c, T.~Szab\'o,
Positional Games,
Birkh\" auser/Springer, Basel, 2014.

\bibitem{henning-2015}
M.A.~Henning, S.~Klav\v{z}ar, D.F.~Rall,
Total version of the domination game,
Graphs Combin.\ 31 (2015) 1453--1462.

\bibitem{henning-2017}
M.A.~Henning, S.~Klav\v{z}ar, D.F.~Rall,
The $4/5$ upper bound on the game total domination number, 
Combinatorica 37 (2017) 223--251.

\bibitem{irsic-2019}
V.~Ir\v{s}i\v{c}, 
Effect of predomination and vertex removal on the game total domination number of a graph,
Discrete Appl.\ Math.\ 257 (2019) 216--225. 

\bibitem{jiang-2019}
Y.~Jiang, M.~Lu, 
Game total domination for cyclic bipartite graphs,
Discrete Appl.\ Math.\ 265 (2019) 120--127.

\bibitem{portier-2025}
J.~Portier, L.~Versteegen, 
A proof of the $3/4$-conjecture for the total domination game,
SIAM J.\ Discrete Math.\ 39 (2025) 1--18.

\bibitem{rahman-2023}
M.L.~Rahman, T.~Watson, 
$6$-uniform {M}aker-{B}reaker game is {PSPACE}-complete,
Combinatorica 43 (2023) 595--612. 

\bibitem{worawannotai-2024}
C.~Worawannotai, K.~Charoensitthichai, 
$4$-total domination game critical graphs,
Discrete Math.\ Algorithms Appl.\ 16 (2024) Paper 2350061. 

\end{thebibliography}
\end{document}